\documentclass[a4paper,10pt]{article}
\usepackage{graphicx}
\usepackage{subfigure}
\usepackage{epstopdf}
\usepackage{booktabs}
\usepackage{threeparttable}
\usepackage{amssymb, amsmath, amsthm, epsfig}
\numberwithin{equation}{section}  % 使公式编号按章节编号
\numberwithin{figure}{section}    % 使图片编号按章节编号
\numberwithin{table}{section}     % 使表格编号按章节编号

\usepackage{tikz}
\usetikzlibrary{decorations.pathmorphing}

\usepackage{enumitem} %使用 enumerate 环境 + \roman

\usepackage{latexsym}
\usepackage{eepic}
\usepackage{epic}
\usepackage{color}

\usepackage{amsmath,textcomp, amssymb, mathrsfs, bm, amsthm, amsfonts}
\usepackage{algpseudocode}
\allowdisplaybreaks[4]
\usepackage{ifpdf}
\usepackage{siunitx}%
\usepackage{chngcntr}

\usepackage{dsfont}
\usepackage{multirow}
\usepackage{bm}
\usepackage{caption}
\usepackage[colorlinks, linkcolor=blue, citecolor=blue]{hyperref}

\allowdisplaybreaks[4]

\theoremstyle{plain}
\newtheorem{lem}{Lemma}[section]
\newtheorem{thm}[lem]{Theorem}
\newtheorem{cor}[lem]{Corollary}

\theoremstyle{definition}
\newtheorem{defn}{Definition}[section]

\theoremstyle{remark}
\newtheorem{rem}{Remark}[section]

\renewcommand{\thefigure}{\thesection.\arabic{figure}}

\newcommand{\ds}{\displaystyle}

\newcommand{\ms}{\medskip}

\def \e{\ensuremath{\mathrm{e}}}
\def \i{\ensuremath{\mathrm{i}}}
\def \d{\ensuremath{\mathrm{d}}}
\begin{document}
\renewcommand{\figurename}{Figure}
\renewcommand{\thesubfigure}{(\alph{subfigure})}
\renewcommand{\thesubtable}{(\alph{subtable})}
\makeatletter
\renewcommand{\p@subfigure}{\thefigure~}

\makeatother
\title{\large\bf Quantitative Analysis of Second-Harmonic Generation in Plasmonic Nanowires: Symmetry Breaking and a Topological Defect}
\author{
Qi Lei\thanks
{School of Mathematics, Hunan University, Changsha 410082, Hunan Province, China.
Email: leiqi@hnu.edu.cn}
\and
Hongyu Liu\thanks
{Department of Mathematics, City University of Hong Kong, Kowloon, Hong Kong, China.
Email:  hongyu.liuip@gmail.com; hongyliu@cityu.edu.hk}
\and
Zhi-Qiang Miao\thanks
{School of Mathematics and Statistics, Central South University, Changsha 410083, Hunan Province, China, and
Department of Mathematics, City University of Hong Kong, Kowloon, Hong Kong, China.
Email: zhiqmiao@csu.edu.cn; zhiqmiao@cityu.edu.hk}
\and
Guang-Hui Zheng\thanks
{Corresponding author. School of Mathematics, Hunan University, Changsha 410082, Hunan Province, China.
Email: zhenggh2012@hnu.edu.cn; zhgh1980@163.com}
}
\date{}
\maketitle

\begin{center}{\bf ABSTRACT}
\end{center}\smallskip
This paper develops a quantitative mathematical theory of second-harmonic generation (SHG) in plasmonic hollow nanowires, highlighting the distinct roles of geometric symmetry breaking and an internal cavity. We consider a hollow nanowire made of a centrosymmetric material, with a concentric circular cavity and a slightly perturbed outer boundary. Relative to the simply connected solid nanowire considered in previous work~\cite{liu2026quantitative}, the cavity introduces a topological defect by changing the cross-section from a disk into a doubly connected domain, whereas the outer-boundary perturbation breaks geometric symmetry without altering the topology. Local inversion-symmetry breaking at the interfaces induces nonlinear surface polarizations. By combining layer-potential techniques with asymptotic expansions in the boundary perturbation parameter, we derive explicit leading- and first-order formulas for the fundamental and second-harmonic fields. For a circular outer boundary, the leading-order second-harmonic field contains only a quadrupolar component, whose amplitude differs from that of a solid nanowire by ($O(\gamma^2)$) as ($\gamma=r_i/r_e\to 0$) and decays as ($O(r_e-r_i)$) in the nonresonant thin-shell limit. For an outer-boundary perturbation proportional to ($\cos(n\theta)$), ($n\geq 3$), the first-order second-harmonic correction contains only modes of orders ($n-2$), $n$, and ($n+2$), yielding dipole, quadrupole, or higher-order multipole radiation for ($n=3$), ($n=4$), or ($n\geq 5$), respectively, provided the lowest-order modal coefficient is nonzero. We further derive explicit plasmonic resonance conditions at the fundamental and second-harmonic frequencies and estimate the singular orders of single and simultaneous resonances, including possibly higher-order singularities arising from coincident resonance factors. These results show that geometric symmetry determines the multipolar radiation modes, while the cavity radius controls their amplitudes and resonant enhancement.

\smallskip
{\bf keywords}: second-harmonic generation; plasmonic resonance; multipolar radiation; symmetry breaking; topological defect; hollow nanowire; layer-potential techniques; asymptotic analysis

\section{Introduction}
Second-harmonic generation (SHG) is a fundamental nonlinear optical process in which an incident field at frequency \(\omega\) induces a response at twice the frequency. Since the early experimental observations of SHG in quartz crystals and at metal surfaces~\cite{brown1965nonlinear,franken1961generation} , its underlying mechanisms and applications have attracted considerable attention, including studies of structured second-harmonic waves~\cite{saltiel2008spatiotemporal,saltiel2008generation}. In addition, SHG and related nonlinear optical effects have found wide-ranging applications in frequency conversion, materials characterization, optical sensing, and high-resolution imaging ~\cite{goodman2018ultrafast,griffiths2020real,li2018optical,nakayama2007tunable,wu2020high}. In centrosymmetric materials, the bulk second-order nonlinear polarization vanishes within the electric-dipole approximation. At material interfaces, however, inversion symmetry is locally broken, allowing nonlinear surface polarizations even when the overall geometry is inversion symmetric. The pioneering theoretical work of Rudnick and Stern~\cite{rudnick1971second} provided a foundation for describing nonlinear optical responses at metal surfaces. Nonlocal bulk contributions may also arise from spatial variations of the fundamental field. This sensitivity to surfaces and interfaces motivates the study of SHG in plasmonic nanostructures, where geometry and material properties can significantly influence the generated radiation.

The development of nanofabrication techniques has stimulated growing interest in SHG from metallic nanoparticles and nanostructures. Localized surface plasmon resonances can enhance the electric field near metallic interfaces~\cite{ammari2013spectral,ando2016plasmon, blaasten2020localization,deng2021mathematical, deng2022plasmon,lei2025hybridization} and thereby amplify nonlinear optical responses~\cite{anceau2003local,bouhelier2003near, krause2004optical,nahata2003enhanced, pendry1999magnetism}. The multipolar structure of second-harmonic radiation is also closely related to the geometric symmetry of the nanostructure. In particular, inversion-symmetric geometries suppress electric-dipole SHG, while quadrupolar and higher-order contributions may remain. Suitable geometric deformations that break inversion symmetry can activate dipolar radiation. SHG has been investigated in a variety of metallic nanostructures, including metal tips, periodically structured metal films, deformed particles, split-ring resonators and their complementary structures, metal--dielectric multilayers, T- and L-shaped particles, nanodimers, and fishnet structures~\cite{bouhelier2003near, canfield2007local,feth2008second,husu2008chiral, kim2008nonlinear,klein2006second,kujala2007multipole, kujala2008multipolar,larciprete2008second, nahata2003enhanced,nappa2005electric, neacsu2005second,shan2006experimental}. These studies highlight the important role of geometric symmetry in determining the allowed multipolar radiation modes. Analytical studies of SHG scattering from spheres and cylinders were carried out within Mie-type frameworks by Dadap et al.~\cite{dadap1999second} and Valencia et al.~\cite{valencia2004second}, clarifying the role of quadrupolar contributions in centrosymmetric structures. Singla and Moch\'an~\cite{singla2019analytical} subsequently developed a perturbation method for an infinitely long cylindrical nanowire with a slightly deformed boundary and derived explicit formulas describing the influence of geometric deformation on the second-harmonic field. More recently, Liu et al.~\cite{liu2026quantitative} employed a symmetry degree based on group theory to establish a quantitative relation between geometric symmetry and the multipolar structure of the radiated field.

The preceding studies provide a mathematical foundation for understanding the influence of geometric symmetry on second-harmonic generation. Extending the solid-nanowire configurations studied in~\cite{liu2026quantitative,singla2019analytical}, we consider a hollow structure with a concentric circular cavity, which introduces an additional interface and renders the cross-section multiply connected. Relative to the simply connected solid-nanowire cross-section studied in the preceding work, we refer to the introduction of this cavity as a topological defect in the geometric sense. This raises the question of how the interaction between the inner and outer interfaces affects the second-harmonic radiation amplitudes and plasmonic resonance conditions. The resulting topological defect introduces an additional geometric parameter, namely the cavity radius, which provides a means of tuning the nonlinear optical response.

In this paper, we study an isolated, infinitely long plasmonic hollow nanowire made of a centrosymmetric material and embedded in vacuum. The nanowire is invariant along the $\hat{x}_3$-direction, and its cross-section consists of a concentric circular cavity of radius $r_i$ and a slightly perturbed outer boundary of reference radius $r_e$. The circular cavity does not, by itself, break global inversion symmetry. The outer-boundary deformation reduces the rotational symmetry and may also break inversion symmetry, depending on its angular structure. The fundamental and second-harmonic fields satisfy transmission conditions on both interfaces, leading to coupled boundary integral equations. The nonlinear surface polarizations and associated surface charges on both boundaries are also taken into account. Working in the quasistatic regime, we combine layer-potential representations with asymptotic expansions in the boundary perturbation parameter $\varepsilon$ to derive explicit leading- and first-order expressions for the fundamental and second-harmonic fields.

Our analysis reveals the distinct roles of the cavity and the outer-boundary perturbation in determining second-harmonic radiation. When the outer boundary is perfectly circular, the exterior leading-order second-harmonic field contains only a quadrupolar component, with no dipole contribution. Writing $\gamma=r_i/r_e$, we show that the quadrupolar amplitude differs from that of the corresponding solid nanowire by $O(\gamma^2)$ as $\gamma\to 0$, whereas it decays as $O(r_e-r_i)$ in the nonresonant thin-shell limit $r_i\to r_e^-$, under the corresponding assumptions on the material parameters. Moreover, the quadrupolar amplitude may depend nonmonotonically on the cavity size. For an outer-boundary perturbation of the form $f(\theta)=r_e\cos(n\theta)$, $n\geq 3$, the first-order second-harmonic correction contains only the angular modes of orders $n-2$, $n$, and $n+2$. Provided that the lowest-order modal coefficient does not vanish, this correction exhibits dipole radiation for $n=3$, quadrupole radiation for $n=4$, and higher-order multipole radiation for $n\geq 5$. These classifications concern the first-order correction, whereas the full second-harmonic field generally also contains the leading-order quadrupolar contribution.

We further derive explicit plasmonic resonance conditions at the fundamental and second-harmonic frequencies and estimate the associated singular orders for both single and simultaneous resonances. These conditions depend explicitly on the geometric ratio $\gamma$, showing that the cavity radius influences the resonance conditions and the enhancement of individual multipolar amplitudes. In particular, coincident resonance factors can lead to higher-order singular behavior, while simultaneous resonances at the two frequencies exhibit different combinations of singular factors. Thus, the cavity radius provides an additional geometric parameter for tuning the multipolar amplitudes and resonance behavior without altering the modal orders determined by the outer-boundary symmetry.

The remainder of the paper is organized as follows. Section~\ref{sec:mathematical setting and main results} introduces the mathematical model and states the main results. Section~\ref{sec:auxiliary results} recalls the necessary tools from group theory and layer-potential theory and establishes well-posedness. Section~\ref{sec:asymptotic analysis of the perturbed fields} develops asymptotic expansions and layer-potential representations for the fundamental and second-harmonic fields. Section~\ref{sec:proofs of the main results} contains the proofs of the main results on multipolar radiation and plasmonic resonances. Finally, Section~\ref{sec:conclusion} concludes the paper and discusses possible extensions.

\section{Mathematical setting and main results}
\label{sec:mathematical setting and main results}
\begin{figure}[htbp]
    \centering

    %==================================================
    % (a) 三维空心纳米线示意图
    %==================================================
    \begin{tikzpicture}[baseline=(current bounding box.center), line join=round, line cap=round]

        %----------------------------
        % 基本参数
        %----------------------------
        \def\R{1.0}         % 外边界参考半径
        \def\a{0.20}        % 外边界扰动幅度
        \def\k{5}           % 扰动模态，对应 cos(5\theta)
        \def\h{5.0}         % 圆柱高度
        \def\vy{0.42}       % 顶/底截面纵向压缩比例（产生透视感）
        \def\Rinner{0.36}   % 内部圆孔半径

        %----------------------------
        % 外边界在极大值处的半径
        % cos(5\theta) 的五个极大值位置：
        % 1: 0°
        % 2: 72°
        % 3: 144°
        % 4: 216°
        % 5: 288°
        %
        % 本次要求：只保留位置 1、3 的实线，
        % 去掉位置 5 的实线
        %----------------------------
        \pgfmathsetmacro{\rmax}{\R+\a}

        % 位置 1：theta = 0°
        \def\thA{0}
        \pgfmathsetmacro{\xoutA}{\rmax*cos(\thA)}
        \pgfmathsetmacro{\youtA}{\rmax*sin(\thA)*\vy}

        % 位置 3：theta = 144°
        \def\thB{144}
        \pgfmathsetmacro{\xoutB}{\rmax*cos(\thB)}
        \pgfmathsetmacro{\youtB}{\rmax*sin(\thB)*\vy}

        %----------------------------
        % 先画顶部外边界（五叶形）
        %----------------------------
        \draw[
            fill=gray!25,
            domain=0:360,
            samples=320,
            smooth,
            variable=\t
        ]
        plot (
            {(\R+\a*cos(\k*\t))*cos(\t)},
            {(\R+\a*cos(\k*\t))*sin(\t)*\vy+\h}
        );

        %----------------------------
        % 再画底部外边界（五叶形）
        %----------------------------
        \draw[
            fill=gray!25,
            domain=0:360,
            samples=320,
            smooth,
            variable=\t
        ]
        plot (
            {(\R+\a*cos(\k*\t))*cos(\t)},
            {(\R+\a*cos(\k*\t))*sin(\t)*\vy}
        );

        %----------------------------
        % 顶部内圆孔（白色填充）
        %----------------------------
        \draw[
            fill=white,
            domain=0:360,
            samples=180,
            smooth,
            variable=\t
        ]
        plot (
            {\Rinner*cos(\t)},
            {\Rinner*sin(\t)*\vy+\h}
        );

        %----------------------------
        % 底部内圆孔（白色填充）
        %----------------------------
        \draw[
            fill=white,
            domain=0:360,
            samples=180,
            smooth,
            variable=\t
        ]
        plot (
            {\Rinner*cos(\t)},
            {\Rinner*sin(\t)*\vy}
        );

        %----------------------------
        % 顶部内边界轮廓
        %----------------------------
        \draw[
            thick,
            domain=0:360,
            samples=180,
            smooth,
            variable=\t
        ]
        plot (
            {\Rinner*cos(\t)},
            {\Rinner*sin(\t)*\vy+\h}
        );

        %----------------------------
        % 底部内边界轮廓
        %----------------------------
        \draw[
            thick,
            domain=0:360,
            samples=180,
            smooth,
            variable=\t
        ]
        plot (
            {\Rinner*cos(\t)},
            {\Rinner*sin(\t)*\vy}
        );

        %==================================================
        % 下面这些线放到最后画，避免被灰色区域遮挡
        %==================================================

        %----------------------------
        % 外边界实线：位置 1
        %----------------------------
        \draw[thick]
        (\xoutA,\youtA)
        --
        (\xoutA,\h+\youtA);

        %----------------------------
        % 外边界实线：位置 3
        %----------------------------
        \draw[thick]
        (\xoutB,\youtB)
        --
        (\xoutB,\h+\youtB);

        %----------------------------
        % 内边界虚线：0° 位置（即右端点）
        % 对应 x = Rinner
        % 顶部和底部交点的 y 坐标都为 0
        %----------------------------
        \draw[dashed]
        (\Rinner,0)
        --
        (\Rinner,\h);

        %----------------------------
        % 内边界虚线：180° 位置（即左端点）
        % 对应 x = -Rinner
        %----------------------------
        \draw[dashed]
        (-\Rinner,0)
        --
        (-\Rinner,\h);

        %----------------------------
        % 坐标轴
        %----------------------------
        \draw[-latex,thick]
        (0,0) -- (1.55,0);
        \node at (1.80,0) {$x_1$};

        \draw[-latex,thick]
        (0,0) -- (0,1.15);
        \node at (0,1.38) {$x_3$};

        \draw[-latex,thick]
        (0,0) -- (-0.68,-0.95);
        \node at (-0.82,-1.14) {$x_2$};

        %----------------------------
        % 入射基频电场
        %----------------------------
        \draw[-latex,thick]
        (-3.10,2.55) -- (-1.55,2.55);

        \draw[thick] (-2.36,2.75) -- (-2.36,2.35);
        \draw[thick] (-2.26,2.75) -- (-2.26,2.35);
        \draw[thick] (-2.16,2.75) -- (-2.16,2.35);

        \node at (-2.25,3.02) {$\bm E_{\omega}^{\mathrm{in}}$};

        %----------------------------
        % 基频散射场
        %----------------------------
        \draw[
            thick,
            domain=0:2.05,
            smooth,
            variable=\x,
            rotate=28,
            shift={(3.20,2.05)}
        ]
        plot ({\x},{0.18*sin(8*\x r)});

        \node at (3.10,3.35) {$\bm E_{\omega}^{\mathrm{sc}}$};

        %----------------------------
        % 二次谐波散射场
        %----------------------------
        \draw[
            thick,
            green!70!black,
            domain=0:2.05,
            smooth,
            variable=\x,
            rotate=-28,
            shift={(0.55,2.00)}
        ]
        plot ({\x},{0.18*cos(8*\x r)});

        \node at (3.00,1.90) {$\bm E_{2\omega}^{\mathrm{sc}}$};

        %----------------------------
        % 子图编号
        %----------------------------
        \node at (0,-1.50) {(a)};

    \end{tikzpicture}
    \qquad
    %==================================================
    % (b) 二维约化模型示意图
    %==================================================
    \begin{tikzpicture}[baseline=(current bounding box.center), scale=1.12, line join=round, line cap=round]

        \def\Rin{0.50}   % 内圆半径

        %----------------------------
        % 外边界：五叶扰动边界
        %----------------------------
        \draw[
            fill=gray!25,
            domain=0:360,
            samples=320,
            smooth,
            variable=\x
        ]
        plot ({\x}:{1+0.2*cos(5*\x)});

        %----------------------------
        % 内部圆孔
        %----------------------------
        \draw[
            fill=white,
            domain=0:360,
            samples=200,
            smooth,
            variable=\x
        ]
        plot ({\x}:{\Rin});

        % 内边界轮廓
        \draw[
            thick,
            domain=0:360,
            samples=200,
            smooth,
            variable=\x
        ]
        plot ({\x}:{\Rin});

        %----------------------------
        % 坐标轴
        %----------------------------
        \draw[-latex,thick]
        (0,0) -- (1.55,0);
        \node at (1.78,0) {$x_1$};

        \draw[-latex,thick]
        (0,0) -- (0,1.55);
        \node at (0,1.78) {$x_2$};

        %----------------------------
        % 区域和边界标注
        %----------------------------
        \node at (1.05,0.50) {$\partial\Omega_\epsilon$};
        \node at (-0.28,0.66) {$\partial D$};
        \node at (0.18,-0.26) {$D$};
        \node[font=\scriptsize] at (-0.64,-0.52) {$\Omega_\epsilon\setminus\overline{D}$};

        %----------------------------
        % 外加电场
        %----------------------------
        \draw[-latex,thick]
        (-2.55,0) -- (-1.05,0);
        \node at (-1.80,0.30) {$\bm E_{\mathrm{ext}}$};

        %----------------------------
        % 基频散射场
        %----------------------------
        \draw[
            thick,
            domain=0:2.05,
            smooth,
            variable=\x,
            rotate=28,
            shift={(1.60,0.05)}
        ]
        plot ({\x},{0.18*sin(8*\x r)});

        \node at (2.68,1.08) {$\bm E_{\omega}^{\mathrm{sc}}$};

        %----------------------------
        % 二次谐波散射场
        %----------------------------
        \draw[
            thick,
            green!70!black,
            domain=0:2.05,
            smooth,
            variable=\x,
            rotate=-28,
            shift={(1.65,-0.02)}
        ]
        plot ({\x},{0.18*cos(8*\x r)});

        \node at (2.64,-1.05) {$\bm E_{2\omega}^{\mathrm{sc}}$};

        %----------------------------
        % 外法向量
        %----------------------------
        \draw[-latex,thick]
        (-0.58,-0.90) -- (-0.85,-1.25);
        \node at (-0.94,-1.38) {$\nu$};

        %----------------------------
        % 子图编号
        %----------------------------
        \node at (0,-2.08) {(b)};

    \end{tikzpicture}

    \caption{
        Schematic illustration of second harmonic generation in a hollow nanowire.
        (a) Three-dimensional model: scattering of a monochromatic time-harmonic incident wave by a hollow nanowire with a perturbed outer boundary.
        (b) The reduced two-dimensional problem, where the domain \(\Omega_\epsilon\setminus\overline{D}\) is an annular region with a perturbed outer boundary and a circular inner boundary.
    }
    \label{fig:schematic}
\end{figure}
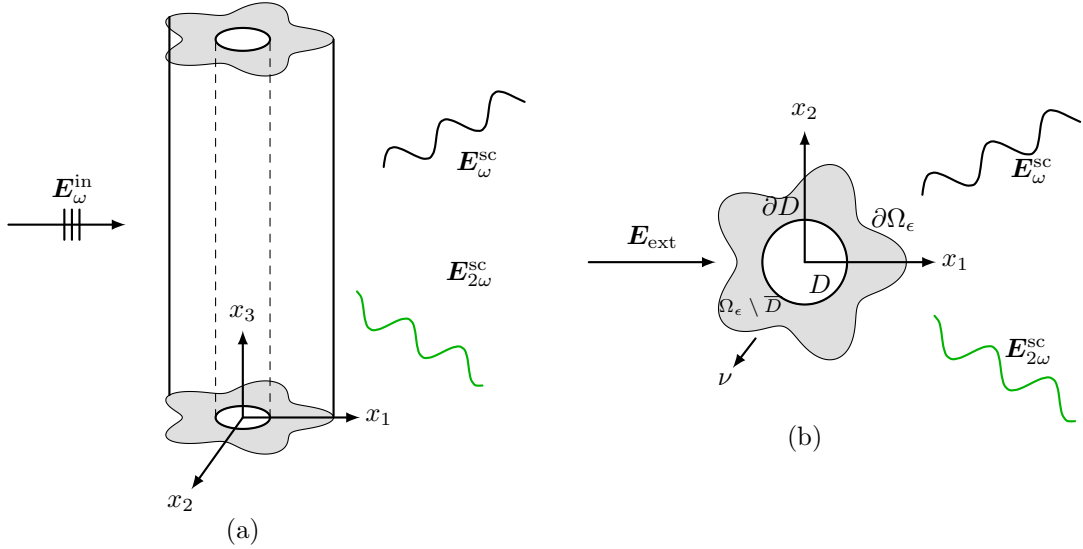

We begin by considering a hollow nanowire composed of a centrosymmetric material, as illustrated in Fig.~\ref{fig:schematic}(a). The nanowire is illuminated by a plane monochromatic laser beam \( \bm E_{\omega}^{\text{in}} \) of frequency \( \omega \), with its polarization lying in the plane of the cross-section. The nanowire is modeled as an isolated, infinitely long cylinder placed in vacuum, exhibiting translational symmetry along its axis (the \( \hat{x}_3 \) direction). The cross-sectional geometry consists of a concentric circular cavity of radius $r_i$ and a slightly perturbed outer boundary. The cavity, viewed as a topological defect relative to the solid-disk cross-section considered in the preceding work, renders the cross-section doubly connected without itself breaking inversion symmetry. By contrast, the outer-boundary deformation reduces rotational symmetry and may also break inversion symmetry, depending on its angular structure. The outer boundary is characterized by a perturbed radius \( \tilde{r}_e \), and the characteristic length scale \( r_e \) (on which the deformation is applied) is assumed to be much smaller than the wavelength of the incident electromagnetic field, allowing us to neglect retardation effects. Under this quasi-static approximation, the plane monochromatic laser beam is reduced to an external electric field \( \bm E_{\text{ext}} \) oscillating in the \( x_1 \)-\( x_2 \) plane, as shown in Fig.~\ref{fig:schematic}(b). Consequently, the problem is treated within the nonretarded (quasi-static) regime.

In this paper, we consider a hollow nanowire made of a centrosymmetric material. Due to the local loss of inversion symmetry within a thin selvedge region near the surface, a nonlinear polarization is induced on both the inner and outer boundaries of the nanowire. This surface nonlinear polarization \( P^s_i \) is expressed as
\begin{align}
    \label{eq:surface-nonlinear-polarization}
    P^s_i = \sum_{jk} \chi_{ijk}^s F_j F_k,
\end{align}
where \( \chi_{ijk}^s \) are the components of the local nonlinear surface susceptibility tensor, characterizing the second-order nonlinear response of the material within the thin selvedge region. The field \( \bm F  \) is constructed using quantities that are continuous across the interface, thereby avoiding ambiguity regarding the position within the selvedge at which the fields are evaluated. Specifically, \( \bm F \) consists of the normal projection of the displacement field and the parallel projection of the electric field, both evaluated at the surface. However, the precise definition of \( \bm F \) depends on the interface under consideration. For the outer boundary \( \partial \Omega_\epsilon \), the nonlinear polarization sheet is located just outside the nanowire, while the fundamental electric field is evaluated just inside the interface. Hence, \( \bm F \) is defined at \( \partial \Omega_\epsilon \) (i.e., just outside the outer boundary) as
\begin{align}
    \label{eq:F-r-e}
    \bm F = \bm E_\omega|_{+} = \varepsilon_\omega \bm E_\omega^{\perp}|_{-} + \bm E_\omega^{\parallel }|_{-},
\end{align}
where the symbols \(\perp\) and \(\parallel\) represent the components normal and tangential to the boundary, respectively. The subscripts $|_{+}$ and $|_{-}$ indicate that the nonlinear polarization sheet is located just above the metal, while the fundamental electric field is evaluated just below the interface. For the inner boundary \( \partial D\), the nonlinear polarization sheet is located just inside the nanowire, while the fundamental electric field is evaluated just outside the interface. Consequently, \( \bm F \) is defined at \( \partial D \) (i.e., just inside the inner boundary) as
\begin{align}
    \label{eq:F-r-i}
    \bm F = \bm E_\omega|_{-} = \varepsilon_\omega \bm E_\omega^{\perp}|_{+} + \bm E_\omega^{ \parallel }|_{+},
\end{align}
%where \( |_{r_i^+} \) denotes evaluation at \( r_i^+ \) (i.e., just outside the inner boundary).
This separate definition ensures that on each interface, the field \( \bm F \) is constructed from quantities that are continuous across that interface, thereby providing a consistent and unambiguous description of the surface nonlinear polarization.  The fundamental electric field \( \bm E_\omega \) can be expressed in terms of the electrostatic potential \( u_\omega \) as \( \bm E_\omega = -\nabla u_\omega \).

We assume that the thickness of the selvedge region is much smaller than the radius of the cylinder, and thus that the surface can be considered as locally flat. We further assume local invariance under rotations around the surface normal. Unless the surface contains structural features with intrinsic chirality, the metal-background interface possesses an isotropic mirror-symmetry plane perpendicular to the surface. Under these circumstances, the surface nonlinear susceptibility \( \chi_{ijk}^s \) has only three independent components, namely, \( \chi_{\perp\perp\perp}^s \), \( \chi^s_{\perp\parallel\parallel} \), and \( \chi_{\parallel\parallel\perp}^s = \chi_{\parallel\perp\parallel}^s \).
The nonlinear polarization induced on the surface of the cylinder is obtained by using \eqref{eq:surface-nonlinear-polarization}--\eqref{eq:F-r-i}. Specifically, its normal component and tangential component are
\begin{align*}
    %\label{eq:Ps-pre}
    P^s_{\perp} & = \chi_{\perp\perp\perp} ^s F_{\perp} F_{\perp}, \\[8pt]
    %\label{eq:Ps-par}
    P^s_{\parallel} & =  2 \chi_{\parallel\perp\parallel} ^s F_{\perp} F_{\parallel}.
\end{align*}
For convenience, we denote \( \chi^s_{\perp} := \chi^s_{\perp\perp\perp} \) and \( \chi^s_{\parallel} := \chi^s_{\parallel\perp\parallel} \) in subsequent sections. For the outer boundary \( \partial\Omega_\epsilon \), the normal and tangential components of the surface nonlinear polarization on the outer boundary, denoted \( P^s_{\perp,e} \) and \( P^s_{\parallel,e} \), can be expressed as
\begin{align}
    \label{eq:Pse-perp}
    P^s_{\perp,e} = \chi^s_{\perp} \varepsilon_\omega^2 \left( \frac{\partial u_\omega}{\partial \nu} \Big|_{-} \right)^2  & \quad \mbox{on } \partial \Omega_\epsilon, \\[4pt]
    \label{eq:Pse-parallel}
    P^s_{\parallel,e} = 2 \chi^s_{\parallel} \varepsilon_\omega \left( \frac{\partial u_\omega}{\partial \nu} \Big|_{-}\right)  \left( \frac{\partial u_\omega}{\partial \tau} \Big|_{-} \right) &\quad \mbox{on } \partial \Omega_\epsilon,
\end{align}
where \( \nu \) and \( \tau \) denote the outward unit normal and tangential vectors to the surface, respectively. For the inner boundary \( \partial D \), the normal and tangential components of the surface nonlinear polarization on the inner boundary, denoted \( P^s_{\perp,i} \) and \( P^s_{\parallel,i} \), are given by
\begin{align}
    \label{eq:Psi-perp}
    P^s_{\perp,i} = \chi^s_{\perp}  \varepsilon_\omega^2 \left( \frac{\partial u_\omega }{\partial \nu} \Big|_{+} \right)^2 & \quad \mbox{on } \partial D, \\[4pt]
    \label{eq:Psi-parallel}
    P^s_{\parallel,i} = 2 \chi^s_{\parallel} \varepsilon_\omega \left( \frac{\partial u_\omega}{\partial \nu} \Big|_{+} \right)  \left( \frac{\partial u_\omega}{\partial \tau}\Big|_{+} \right) & \quad \mbox{on } \partial D.
\end{align}
The variation of the tangential component of the nonlinear surface polarization along the surface yields another contribution to the surface charge $\sigma^s$, where we use the superscript $s$ to denote its surface origin. For the outer boundary $\partial\Omega_\epsilon$, the divergence of the tangential nonlinear surface polarization $\bm P^s_{\parallel,e}$ along the outer surface contributes to the surface charge $\sigma^s_e$, which is defined as
\begin{align}
    \label{eq:sigma-s-e}
    \sigma^s_e = - \nabla_{\parallel} \cdot \bm P^s_{\parallel,e},
\end{align}
where $\bm P^s_{\parallel,e} = P^s_{\parallel,e} \tau_e$, $\tau_e$ is the unit tangential vector on the outer boundary, $\nabla_{\parallel}$ denotes the surface gradient operator on the outer surface. For the inner boundary $\partial D$, the divergence of the tangential nonlinear surface polarization $\bm P^s_{\parallel,i}$ along the inner surface contributes to the surface charge $\sigma^s_i$, which is defined as
\vspace{6pt}
\begin{equation}
    \label{eq:sigma-s-i}
    \sigma^s_i = - \nabla_{\parallel} \cdot \bm P^s_{\parallel,i},
\vspace{6pt}
\end{equation}
where $\bm P^s_{\parallel,i} = P^s_{\parallel,i} \tau_i$, $\tau_i$ is the unit tangential vector on the inner boundary.

To mathematically state the problem, let \( D \) and \( \Omega \) be two concentric disks in \( \mathbb{R}^2 \), centered at the origin, with radii \(0<r_i<r_e\). We define the bounded domain \( \Omega_\epsilon \) as a perturbed version of the reference disk \( \Omega \)
%such that \( \Omega_\epsilon \setminus \overline{D} \) is noncentrosymmetric,
with its boundary given by
\begin{align}
    \label{eq:perturbed-boundary}
    \partial \Omega_\epsilon  = \left\{ \tilde{x} \mid \tilde{x} = x + \epsilon \, f(x) \, \nu(x), \quad x \in \partial \Omega \right\},
\end{align}
where \( \epsilon \ll 1 \) is a small perturbation amplitude, the function \( f(x) \) is referred to as the shape function, which governs the geometric shape of the perturbed domain \( \Omega_\epsilon \). In polar coordinates, the perturbed radius is expressed as
\begin{align*}
    \tilde{r}_e = r_e + \epsilon \, f(\theta).
\end{align*}
The permittivities of the nanowire region \( \Omega_\epsilon \setminus \overline{D} \) at the fundamental and second harmonic frequencies are \( \varepsilon_\omega \) and \( \varepsilon_{2\omega} \), respectively, while the background medium, i.e., \(D \cup \left(\mathbb{R}^2 \setminus \overline{\Omega}_\epsilon\right) \) is assigned a permittivity of \( \varepsilon_m = 1 \). Throughout this paper, we assume that $D$ and $\Omega_\epsilon$ are of class $C^{2}$. Let $H(x) = -E x_1$ be the harmonic function in $\mathbb{R}^2$, denoting the background electrostatic potential. The governing equations are modeled as follows
\begin{align}
    \label{eq:governing-equation}
    \renewcommand{\arraystretch}{1.5}
    \begin{array}{l l l}
    \left\{
    \begin{array}{l}
        \Delta u_\omega = 0 \\[8pt]
        u_\omega|_+ = u_\omega|_- \\[8pt]
        u_\omega|_+ = u_\omega|_- \\[8pt]
        \varepsilon_\omega
        \dfrac{\partial u_\omega}{\partial \nu_i}\bigg|_+
        =
        \dfrac{\partial u_\omega}{\partial \nu_i}\bigg|_- \\[8pt]
        \dfrac{\partial u_\omega}{\partial \nu_e}\bigg|_+
        =
        \varepsilon_\omega
        \dfrac{\partial u_\omega}{\partial \nu_e}\bigg|_- \\[8pt]
        u_\omega-H(x)=O(|x|^{-1})
    \end{array}
    \right.
    &
    \begin{array}{l}
        \Delta u_{2\omega}=0 \\[8pt]
        u_{2\omega}|_+-u_{2\omega}|_-
        =
        -4\pi P_{\perp,i}^{s} \\[8pt]
        u_{2\omega}|_+-u_{2\omega}|_-
        =
        4\pi P_{\perp,e}^{s} \\[8pt]
        \varepsilon_{2\omega}
        \dfrac{\partial u_{2\omega}}{\partial \nu_i}\bigg|_+
        -
        \dfrac{\partial u_{2\omega}}{\partial \nu_i}\bigg|_-
        =
        4\pi\sigma_i^{s} \\[8pt]
        \dfrac{\partial u_{2\omega}}{\partial \nu_e}\bigg|_+
        -
        \varepsilon_{2\omega}
        \dfrac{\partial u_{2\omega}}{\partial \nu_e}\bigg|_-
        =
        -4\pi\sigma_e^{s} \\[8pt]
        u_{2\omega}=o(1)
    \end{array}
    &
    \begin{array}{l}
        \textnormal{in }
        \mathbb{R}^2
        \setminus
        \bigl(\partial D \cup \partial \Omega_\epsilon \bigr),
        \\[8pt]
        \textnormal{on }\partial D,
        \\[8pt]
        \textnormal{on }\partial \Omega_\epsilon,
        \\[8pt]
        \textnormal{on }\partial D,
        \\[8pt]
        \textnormal{on }\partial\Omega_\epsilon,
        \\[8pt]
        \textnormal{as }|x|\to+\infty.
    \end{array}
    \end{array}
\end{align}
where \(\nu_i\) and \(\nu_e\) denote the outward unit normal vectors  on \(\partial D\) and  \( \partial \Omega_\epsilon \), respectively, and \(\partial / \partial \nu_i\) and  \( \partial / \partial\nu_e \) denote the corresponding normal derivatives. For \(x\in\partial D\cup\partial\Omega_\epsilon\), let
\(\nu=\nu_i\) on \(\partial D\) and \(\nu=\nu_e\) on
\(\partial\Omega_\epsilon\). We define
\begin{align*}
    u|_{\pm}(x) := \lim_{t\to0^+}u\bigl(x\pm t\nu(x)\bigr),
\end{align*}
and
\begin{align*}
    \left. \frac{ \partial u}{ \partial \nu} \right|_{\pm}(x) := \lim_{t \to 0^+} \left \langle \nabla u \bigl( x \pm t \nu(x) \bigr), \nu(x) \right\rangle,
\end{align*}
provided that the limits exist.

To describe the far-field behavior of the second-harmonic field, we first introduce the following definition~\cite{liu2026quantitative}.

\begin{defn}
Let $u_{2\omega}$ be harmonic in the exterior of a disk and satisfy $u_{2\omega}(x)=o(1)$ as $|x|\to\infty$. Its far-field expansion takes the form
\begin{align*}
    u_{2\omega}(r,\theta) =\sum_{m=1}^{\infty} \frac{1}{r^m} \bigl(a_m\cos(m\theta) + b_m\sin(m\theta)\bigr).
\end{align*}
We say that $u_{2\omega}$ exhibits
$2^m$-multipole radiation if
\begin{align*}
    a_j=b_j=0,\quad 1\leq j<m, \qquad (a_m,b_m)\neq(0,0).
\end{align*}
In particular, $m=1$ and $m=2$ correspond to dipole and quadrupole radiation, respectively.
\end{defn}

This definition will be used below to characterize the far-field decay of the second-harmonic field.

We now present the main results of this paper, whose proofs are given in Section~\ref{sec:proofs of the main results}. We begin by examining the multipolar structure of the second-harmonic field in the absence of outer-boundary perturbations.

\begin{thm}
\label{thm:mainthm1}
Let the hollow nanowire be defined by \eqref{eq:perturbed-boundary}, with the inner boundary \( \partial D\) being a circle of radius \( r_i \). Assume that the fundamental and second-harmonic fields satisfy the governing system \eqref{eq:governing-equation}. If the outer boundary is a perfect circle, namely, \(\epsilon=0\), then the second-harmonic field satisfies
\begin{align*}
    u_{2\omega}(x) = \mathcal{O} \bigl( |x|^{-2} \bigr) \qquad \text{as } |x| \to +\infty.
\end{align*}
Consequently, \(u_{2\omega}\) contains no dipole component, and dipole radiation remains forbidden even in the presence of the concentric circular cavity.
\end{thm}

We next quantify the influence of the cavity size on the leading-order quadrupolar amplitude. For notational convenience, we define $\varepsilon_{\ell\omega\pm1}:=\varepsilon_{\ell\omega}\pm1$ for $\ell=1,2$.

\begin{cor}\label{cor:inner-radius}
Under the assumptions of Theorem~\ref{thm:mainthm1}, the exterior leading-order second-harmonic field admits the representation
\begin{align*}
    u_{2\omega}^{(0)}(x) = \hat{u}(\gamma)\frac{r_e^2}{r^2}\cos(2\theta), \qquad r>r_e,
\end{align*}
The quadrupolar amplitude \(\hat{u}(\gamma)\) has the following properties.

\begin{enumerate}[
    label=\textup{(\roman*)},
    leftmargin=*,
    itemsep=5pt,
    topsep=4pt]

    \item For \(r_i=0\), i.e., \(\gamma=0\),
    \begin{align*}
        \hat{u}(0) = \frac{ 8\pi E^2\varepsilon_\omega \bigl( \chi_\perp^s \varepsilon_\omega \varepsilon_{2\omega} +2\chi_\parallel^s \bigr) }{ \varepsilon_{\omega+1}^2 \varepsilon_{2\omega+1} }.
    \end{align*}

    \item For \(0<\gamma\ll1\),
    \begin{align*}
        \hat{u}(\gamma) = \hat{u}(0) + \mathcal O(\gamma^2).
    \end{align*}
    Thus, a small concentric cavity produces only a second-order correction to the quadrupolar amplitude.

    \item If \(\gamma\to1^-\) and
    \(\varepsilon_\omega \varepsilon_{2\omega} \neq0\), then
    \begin{align*}
        \hat{u}(\gamma) = \hat{u}_*(1-\gamma^2) + \mathcal O\bigl((1-\gamma^2)^2\bigr),
    \end{align*}
    where
    \begin{align*}
        \hat{u}_* = \pi E^2 \left[ \chi_\perp^s (\varepsilon_\omega + \varepsilon_{2\omega}-1) + \frac{\chi_\parallel^s} {\varepsilon_\omega \varepsilon_{2\omega}} \bigl( \varepsilon_\omega^2 \varepsilon_{2\omega} -2\varepsilon_\omega \varepsilon_{2\omega} +2\varepsilon_\omega + \varepsilon_{2\omega} \bigr) \right].
    \end{align*}
    Hence
    \begin{align*}
        \hat{u}(\gamma) = \mathcal O(r_e-r_i) \qquad \text{as }r_i \to r_e^-.
    \end{align*}

    \item For \(0<\gamma<1\), $\hat{u}(\gamma)$ may depend nonmonotonically on $\gamma$, and its variation may be strongly affected by the plasmonic resonances described in Theorem~\ref{thm:mainthm3}.
\end{enumerate}
\end{cor}

We now consider the effect of outer-boundary perturbations on the multipolar structure of the second-harmonic field.

\begin{thm}
\label{thm:mainthm2}
Under the geometric and material setting described in~\eqref{eq:perturbed-boundary}--\eqref{eq:governing-equation} , suppose that the outer-boundary perturbation is given by the shape function
\begin{align*}
    f(\theta)=r_e\cos(n\theta), \qquad n\geq3.
\end{align*}
Then the first-order correction \(u_{2\omega}^{(1)}\) to the second-harmonic field in the exterior region contains only the multipole modes of orders \(n-2\), \(n\), and \(n+2\). Provided that the coefficient of the lowest-order mode does not vanish, its leading far-field behavior is characterized as follows:
\begin{itemize}
    \item If \(n=3\), the leading term decays as
    \(\mathcal{O}(|x|^{-1})\), and hence dipole radiation is excited.

    \item If \(n=4\), the leading term decays as
    \(\mathcal{O}(|x|^{-2})\), corresponding to quadrupole radiation.

    \item If \(n\geq5\), the leading term decays as
    \(\mathcal{O}(|x|^{-(n-2)})\), corresponding to a
    \(2^{\,n-2}\)-pole radiation field.
\end{itemize}
Here, \( u_{2\omega}^{(1)} \) denotes the first-order correction in the asymptotic expansion of the second-harmonic field \( u_{2\omega} \) given in \eqref{eq:u-2omega-expansion}.
\end{thm}

\begin{rem}
These classifications concern the first-order correction \(u_{2\omega}^{(1)}\); the full second-harmonic field generally also contains the leading-order quadrupolar contribution.
\end{rem}

We further investigate the plasmonic resonance behavior of the first-order second-harmonic field.
\begin{thm}
\label{thm:mainthm3}
Under the assumptions of Theorem~\ref{thm:mainthm2}, fix $n\geq3$ and $0<\gamma<1$, the exterior first-order second-harmonic field has the following resonance orders:
\begin{enumerate}[label=\textup{(\roman*)},leftmargin=*,itemsep=8pt,topsep=5pt]
    \item When only one resonance condition is approached, the resonance orders are as follows.
    \begin{itemize}[leftmargin=2em,itemsep=6pt]
        \item If $\varepsilon_{\omega+1} / \varepsilon_{\omega-1} \to \pm \gamma$, $\pm\gamma^{n-1}$, or $\pm\gamma^{n+1}$, the corresponding resonance orders are
        \begin{align*}
        \mathcal O\!\left( \left|\varepsilon_{\omega+1} \mp \gamma \varepsilon_{\omega-1} \right|^{-2}\right), \qquad \mathcal O\! \left( \left| \varepsilon_{\omega+1} \mp \gamma^{n-1} \varepsilon_{\omega-1} \right|^{-1} \right),\qquad \mathcal O\!\left( \left|\varepsilon_{\omega+1} \mp \gamma^{n+1} \varepsilon_{\omega-1} \right|^{-1}\right).
        \end{align*}

        \item If $\varepsilon_{2\omega+1} / \varepsilon_{2\omega-1} \to \pm \gamma^2$, the resonance order is
        \begin{align*}
        \mathcal O\! \left( \left| \varepsilon_{2\omega+1} \mp \gamma^{2} \varepsilon_{2\omega-1} \right|^{-(1+\delta_{n=4})} \right).
        \end{align*}
        If $\varepsilon_{2\omega+1} / \varepsilon_{2\omega-1} \to \pm \gamma^{n+p_0}$, where $p_0\in\{-2,0,2\}$ and $n+p_0\neq2$, the resonance order is
        \begin{align*}
        \mathcal O\! \left( \left| \varepsilon_{2\omega+1} \mp \gamma^{n+p_0} \varepsilon_{2\omega-1} \right|^{-1} \right).
        \end{align*}
    \end{itemize}

        \item When one fundamental-frequency and one second-harmonic resonance condition are approached simultaneously, the resonance orders are as follows.
    \begin{itemize}[leftmargin=2em,itemsep=6pt]
        \item Suppose that $\varepsilon_{\omega+1}/\varepsilon_{\omega-1}\to\pm\gamma$. If $\varepsilon_{2\omega+1}/\varepsilon_{2\omega-1}\to\pm\gamma^2$, the resonance order is
        \begin{align*}
        \mathcal O\! \left( \left| \varepsilon_{\omega+1} \mp \gamma \varepsilon_{\omega-1} \right|^{-2} \left| \varepsilon_{2\omega+1} \mp \gamma^{2} \varepsilon_{2\omega-1} \right|^{-(1+\delta_{n=4})}\right).
        \end{align*}
        If $\varepsilon_{2\omega+1}/\varepsilon_{2\omega-1}\to\pm\gamma^{n+p_0}$, where $p_0\in\{-2,0,2\}$ and $n+p_0\neq2$, the resonance order is
        \begin{align*}
        \mathcal O\! \left( \left| \varepsilon_{\omega+1} \mp \gamma \varepsilon_{\omega-1} \right|^{-2} \left| \varepsilon_{2\omega+1} \mp \gamma^{n+p_0} \varepsilon_{2\omega-1} \right|^{-1}\right).
        \end{align*}

        \item Suppose that $\varepsilon_{\omega+1} / \varepsilon_{\omega-1} \to \pm \gamma^{n-1}$. If $\varepsilon_{2\omega+1} / \varepsilon_{2\omega-1} \to \pm \gamma^2$, the resonance order is
        \begin{align*}
        \mathcal O\! \left( \big( \left| \varepsilon_{\omega+1} \mp \gamma^{n-1} \varepsilon_{\omega-1} \right|^{-1} + \left| \varepsilon_{2\omega+1} \mp \gamma^{2} \varepsilon_{2\omega-1} \right|^{-1} \big) \left| \varepsilon_{2\omega+1} \mp \gamma^{2} \varepsilon_{2\omega-1} \right|^{-\delta_{n=4}} \right).
        \end{align*}
        If $\varepsilon_{2\omega+1} / \varepsilon_{2\omega-1} \to \pm \gamma^{n+p_0}$ and $n+p_0\neq2$, the resonance order is
        \begin{align*}
         \begin{cases}
         \mathcal O\! \left( \left| \varepsilon_{\omega+1} \mp \gamma^{n-1} \varepsilon_{\omega-1} \right|^{-1} \left| \varepsilon_{2\omega+1} \mp \gamma^{n+p_0} \varepsilon_{2\omega-1} \right|^{-1}\right), &p_0\in\{-2,0\},\\[6pt]
         \mathcal O\! \left( \left| \varepsilon_{\omega+1} \mp \gamma^{n-1} \varepsilon_{\omega-1} \right|^{-1} + \left| \varepsilon_{2\omega+1} \mp \gamma^{n+p_0} \varepsilon_{2\omega-1} \right|^{-1}\right),&p_0=2.
         \end{cases}
        \end{align*}

        \item Suppose that $\varepsilon_{\omega+1} / \varepsilon_{\omega-1} \to \pm \gamma^{n+1}$. If $\varepsilon_{2\omega+1} / \varepsilon_{2\omega-1} \to \pm \gamma^2$, the resonance order is
        \begin{align*}
        \mathcal O\! \left( \left| \varepsilon_{\omega+1} \mp \gamma^{n+1} \varepsilon_{\omega-1} \right|^{-1} + \left| \varepsilon_{2\omega+1} \mp \gamma^{2} \varepsilon_{2\omega-1} \right|^{-(1+\delta_{n=4})} \right).
        \end{align*}
        If $\varepsilon_{2\omega+1} / \varepsilon_{2\omega-1} \to \pm \gamma^{n+p_0}$ and $n+p_0\neq2$, the resonance order is
        \begin{align*}
         \begin{cases}
         \mathcal O\! \left( \left| \varepsilon_{\omega+1} \mp\gamma^{n+1} \varepsilon_{\omega-1} \right|^{-1} + \left| \varepsilon_{2\omega+1} \mp \gamma^{n-2} \varepsilon_{2\omega-1} \right|^{-1} \right),&p_0=-2,\\[6pt]
         \mathcal O\! \left( \left| \varepsilon_{\omega+1} \mp \gamma^{n+1} \varepsilon_{\omega-1} \right|^{-1} \left| \varepsilon_{2\omega+1} \mp \gamma^{n+p_0} \varepsilon_{2\omega-1} \right|^{-1} \right), &p_0\in\{0,2\}.
         \end{cases}
        \end{align*}
    \end{itemize}
\end{enumerate}
Here $\delta_{n=4}=1$ if $n=4$ and $0$ otherwise.
\end{thm}

\begin{rem}
Theorems~\ref{thm:mainthm2} and~\ref{thm:mainthm3} describe two complementary aspects of the first-order second-harmonic field. The symmetry index $n$ determines the allowed multipolar modes, whereas plasmonic resonances influence their amplitudes. Provided that the lowest-order modal coefficient is nonzero, the leading radiation is dipolar for $n=3$, quadrupolar for $n=4$, and of higher multipole order for $n\geq 5$. Resonances may enhance these contributions but do not alter the symmetry-based selection of angular modes. These classifications concern the first-order correction, while the full second-harmonic field generally also contains the leading-order quadrupolar contribution.
\end{rem}

\begin{rem}
Theorem~\ref{thm:mainthm3} reveals a hierarchy of possible singular orders. At the fundamental frequency, the condition $\varepsilon_{\omega+1}/\varepsilon_{\omega-1} =\pm\gamma$ corresponds to a second-order singular factor, whereas the other fundamental-frequency resonance conditions involve first-order factors. When resonance occurs only at the second-harmonic frequency, the singular order is generally at most one; when $n=4$ and the source-mode condition $\varepsilon_{2\omega+1}/\varepsilon_{2\omega-1}\to\pm\gamma^2$ is approached, the source and lowest output denominators coincide, giving an upper bound of order two. For simultaneous resonances, the singular factors may combine multiplicatively or additively, depending on the corresponding modal contributions. In particular, coincident resonance conditions may produce higher-order singular behavior. The resulting estimates provide upper bounds on the possible enhancement, which may be reduced by vanishing coefficients or cancellations.
\end{rem}

\begin{rem}
The geometry-dependent resonance conditions in Theorem~\ref{thm:mainthm3} involve the radius ratio $\gamma$. At the fundamental frequency, these conditions depend on $\gamma$, $\gamma^{n-1}$ and $\gamma^{n+1}$, whereas at the second-harmonic frequency they involve $\gamma^2$, $\gamma^{n-2}$, $\gamma^n$, and $\gamma^{n+2}$. In particular, when $n=4$, the identity $\gamma^{n-2}=\gamma^2$ causes the lowest output-mode denominator to coincide with the source-mode denominator, potentially producing higher-order singular behavior. Thus, varying the cavity radius provides a means of tuning the geometry-dependent resonance conditions and may facilitate simultaneous resonances for suitable material parameters.
\end{rem}

\begin{rem}
For the outer-boundary perturbation $f(\theta) = r_e \cos(n\theta)$, we have $f(\theta+\pi) = (-1)^n f(\theta)$. Thus, for a nonzero perturbation amplitude, the cross-section preserves inversion symmetry when $n$ is even and breaks it when $n$ is odd. However, inversion-symmetry breaking alone does not guarantee dipole radiation. Indeed, Theorem~\ref{thm:mainthm2} shows that a dipole component can occur in the first-order second-harmonic field only for $n=3$. Hence, for odd $n\geq 5$, inversion symmetry is broken, but no dipole mode is generated at first order.
\end{rem}

\section{Auxiliary results}
\label{sec:auxiliary results}
In this section, we collect several auxiliary results from group theory and layer-potential theory that will be used in the subsequent analysis.

\subsection{Group theory and symmetry degree}
\label{sec:group-theory}
In this subsection, we recall some basic notions from group theory that are essential for quantifying the symmetry of the nanowire cross-section. A detailed exposition can be found in~\cite{liu2026quantitative}.

Let $D_n$ denote the dihedral group of order $2n$, consisting of $n$ rotations and $n$ reflections that map a regular $n$-gon onto itself. A closed curve $\Gamma$ (or the domain it encloses) is said to be invariant under $D_n$ if $d\Gamma=\Gamma$ for every $d\in D_n$. For a shape invariant under $D_n$ but not under any larger dihedral group, its symmetry degree is defined as $2n$, following~\cite{liu2026quantitative}. We refer to $n$ as the symmetry index, so that the symmetry degree is twice the symmetry index.

In this work, the outer boundary of the hollow nanowire is described in polar coordinates by
\begin{align*}
    \widetilde{r}_e(\theta) =r_e+\varepsilon f(\theta), \qquad f(\theta)=r_e\cos(n\theta),
\end{align*}
with $n\geq 3$. For a nonzero perturbation parameter $\varepsilon$, such a boundary is invariant under the dihedral group $D_n$ and therefore has symmetry degree $2n$. Among the perturbations considered here, $n=3$ corresponds to the smallest symmetry degree, namely $6$, while larger values of $n$ correspond to higher symmetry degrees. This notion will be used to characterize the multipolar radiation of the second-harmonic field.

\subsection{Layer-potential formulation and well-posedness}
\label{sec:layer-potential-formulation-and-well-posedness}
We denote by $G(x, y)$ the fundamental solution of the Laplacian in free space. In $\mathbb{R}^2$, it is given by
\begin{align*}
    G(x, y) = \frac{1}{2\pi} \ln |x - y|.
\end{align*}

Let $D \subset \mathbb{R}^2$ be a bounded domain of class $C^2$. The single-layer potential $\mathcal{S}_D$ and double-layer potential $\mathcal{D}_D$ associated with $D$ are defined, for a density function $\varphi \in L^2(\partial D)$, by
\begin{align*}
    \mathcal{S}_D[\varphi](x) = \int_{\partial D} G(x, y) \varphi(y) \, \d s(y), \quad  x\in \mathbb{R}^2,
\end{align*}
and
\begin{align*}
    \mathcal{D}_D[\varphi](x) = \int_{\partial D} \frac{\partial G(x, y)}{\partial \nu(y)} \varphi(y) \, \d s(y), \quad x\in \mathbb{R}^2\setminus D,
\end{align*}
respectively, where $\nu$ denotes the outward unit normal at $ \partial D$. For $x \in \partial D$, $\mathcal{S}_D[\varphi](x)$ and $\mathcal{D}_D[\varphi](x)$ are referred to as the single-layer potential and the double-layer potential operators, respectively. The following jump relations are well known and will be used extensively:
\begin{align*}
    %\label{jump-relation-S}
    \frac{ \partial \mathcal{S}_{D} [\varphi] }{\partial \nu} \bigg|_{\pm}(x) & = \left( \pm \frac{1}{2} \mathcal{I} + \mathcal{K}_{D}^* \right) [\varphi](x), \quad x \in \partial D, \\
    %\label{jump-relation-D}
    \mathcal{D}_{D} [\varphi]|_{\pm}(x) & = \left( \mp \frac{1}{2} \mathcal{I} + \mathcal{K}_{D} \right)[\varphi](x), \quad x \in \partial D,
\end{align*}
where \( \mathcal{K}_{D} \) is the boundary integral operator defined by
\begin{align*}
    \mathcal{K}_{D}[\varphi](x) = \int_{\partial D} \frac{\partial G(x-y)}{\partial \nu(y)} \varphi(y) \d s(y),
\end{align*}
and \( \mathcal{K}_{D}^* \) is the \( L^2 \)-adjoint of \( \mathcal{K}_{D} \), i.e.,
\begin{align*}
    \mathcal{K}_{D}^* [\varphi](x) = \int_{\partial D} \frac{ \partial G(x-y)}{\partial \nu(x)} \varphi(y) \d s(y).
\end{align*}
Both $\mathcal{K}_D$ and $\mathcal{K}_D^*$ are also known as Neumann--Poincar\'e (NP) operators.

The following lemma provides explicit expressions for the single- and double-layer potentials on a disk \cite{ammari2013spectral,ammari2013mathematical}, which will be essential in the subsequent analysis.
\begin{lem}
\label{Layer-potentials-on-disk}
Let \(D\subset \mathbb{R}^2\) be the disk of radius \(r_0\) centered at the origin. Then, for every nonzero integer \(n\), the following identities hold in polar coordinates \((r,\theta)\).
\begin{enumerate}[label=\normalfont(\roman*)]
    \item The single-layer potential associated with the density \(\e^{\i n\theta}\) is given by
    \begin{equation}
    \label{eq:single-layer-explicit}
    \mathcal{S}_D[\e^{\i n\theta}](x) =
    \begin{cases}
    -\dfrac{r_0}{2|n|}
    \left(\dfrac{r}{r_0}\right)^{|n|}
    \e^{\i n\theta},
    & r<r_0, \\[12pt]
    -\dfrac{r_0}{2|n|}
    \left(\dfrac{r_0}{r}\right)^{|n|}
    \e^{\i n\theta},
    & r>r_0,
    \end{cases}
    \end{equation}
    and its radial derivative satisfies
    \begin{equation}
    \label{eq:single-layer-derivative-explicit}
    \frac{\partial}{\partial r}
    \mathcal{S}_D[\e^{\i n\theta}](x) =
    \begin{cases}
    -\dfrac{1}{2}
    \left(\dfrac{r}{r_0}\right)^{|n|-1}
    \e^{\i n\theta},
    & r<r_0, \\[12pt]
    \dfrac{1}{2}
    \left(\dfrac{r_0}{r}\right)^{|n|+1}
    \e^{\i n\theta},
    & r>r_0.
    \end{cases}
    \end{equation}

    \item The double-layer potential associated with the density \(\e^{\i n\theta}\) is given by
    \begin{equation}
    \label{eq:double-layer-explicit}
    \mathcal{D}_D[\e^{\i n\theta}](x) =
    \begin{cases}
    \dfrac{1}{2}
    \left(\dfrac{r}{r_0}\right)^{|n|}
    \e^{\i n\theta},
    & r<r_0, \\[12pt]
    -\dfrac{1}{2}
    \left(\dfrac{r_0}{r}\right)^{|n|}
    \e^{\i n\theta},
    & r>r_0,
    \end{cases}
    \end{equation}
    while its radial derivative is given by
    \begin{equation}
    \label{eq:double-layer-derivative-explicit}
    \frac{\partial}{\partial r}
    \mathcal{D}_D[\e^{\i n\theta}](x) =
    \begin{cases}
    \dfrac{|n|}{2r}
    \left(\dfrac{r}{r_0}\right)^{|n|}
    \e^{\i n\theta},
    & r<r_0, \\[12pt]
    \dfrac{|n|}{2r}
    \left(\dfrac{r_0}{r}\right)^{|n|}
    \e^{\i n\theta},
    & r>r_0.
    \end{cases}
    \end{equation}

    \item For the constant density \(1\), the double-layer potential satisfies
    \begin{equation}
    \label{eq:double-layer-constant-density-explicit}
    \mathcal{D}_D[1](x) =
    \begin{cases}
    1, & r<r_0, \\[8pt]
    0, & r>r_0.
    \end{cases}
    \end{equation}

    \item The \textnormal{NP} operators \(\mathcal{K}_D^*\) and \(\mathcal{K}_D\) satisfy
    \begin{equation}
    \label{eq:NP-operator-explicit}
    \mathcal{K}_D^*[\e^{\i n\theta}] = \mathcal{K}_D[\e^{\i n\theta}] =
    \begin{cases}
    0, &  n\neq 0, \\[8pt]
    \dfrac{1}{2}, & n= 0.
    \end{cases}
    \end{equation}
\end{enumerate}
\end{lem}

To establish the well-posedness of the coupled system and to derive the representation formulas for the solutions to the governing equations \eqref{eq:governing-equation}, we first recall the following lemma on the spectral properties of the Neumann--Poincar{\'e}-type operator \(\mathbb{K}^*\)~\cite{ammari2013spectral,kang2015layer}.
\begin{lem}
%\label{lem:invertibility-K*}
Let the Neumann--Poincar{\'e}-type operator \(\mathbb{K}^* : \mathcal{H} \to \mathcal{H}\) be defined by
\begin{align*}
    \mathbb{K}^* := \begin{bmatrix}
    \ds -\mathcal{K}_{D}^* & \ds -\dfrac{\partial }{\partial\nu_D} \mathcal{S}_{\Omega_\epsilon} \\[12pt]
    \ds \dfrac{\partial}{\partial\nu_{\Omega_\epsilon}} \mathcal{S}_D & \ds  \mathcal{K}_{\Omega_\epsilon }^* \end{bmatrix},
\end{align*}
where $\mathcal{H} = L^2(\partial D) \times L^2(\partial \Omega_\epsilon)$. Then the spectrum of $\mathbb{K}^*$ satisfies $\sigma(\mathbb{K}^*) \in  [-1/2, 1/2]$. Consequently, for any $\lambda$ with $|\lambda| > 1/2$, the operator $\lambda \mathbb{I} - \mathbb{K}^*$ is invertible on $\mathcal{H}$.
\end{lem}

We now extend above results to the configuration with two disjoint boundaries $\partial D$ and $\partial \Omega_\epsilon$, as described in the following theorem.
\begin{thm}\label{well-posedness}
Let $u_\omega,u_{2\omega} \in C^2\!\left( \mathbb R^2\setminus (\partial D\cup\partial\Omega_\epsilon) \right),$ be the classical solutions of \eqref{eq:governing-equation}. Then $u_\omega$ can be represented as
\begin{align*}
    u_\omega(x) = H(x) + \mathcal{S}_{D} [\varphi_i](x) + \mathcal{S}_{\Omega_\epsilon} [\varphi_e](x), \quad x \in \mathbb{R}^2,
\end{align*}
where the density functions \(\varphi_i \in L_0^2(\partial D)\) and \(\varphi_e \in L_0^2(\partial \Omega_\epsilon)\) satisfy the coupled system
\begin{align*}
    \begin{cases}
    \ds (\lambda_\omega \mathcal{I} + \mathcal{K}_{ D}^*) [\varphi_i] +  \frac{\partial \mathcal{S}_{ \Omega_\epsilon} [\varphi_e] }{\partial\nu_D}= -\dfrac{\partial H }{ \partial \nu_{D}}  & \textnormal{on } \partial D, \ms \\[8pt]
    \ds (\lambda_\omega \mathcal{I} - \mathcal{K}_{ \Omega_\epsilon}^*) [\varphi_e] - \frac{\partial \mathcal{S}_D [\varphi_i] }{ \partial\nu_{ \Omega_\epsilon} }= \dfrac{\partial H}{\partial \nu_{\Omega_\epsilon}}  & \textnormal{on } \partial \Omega_\epsilon,
    \end{cases}
\end{align*}
with $\lambda_\omega$ given by
\begin{align}
    \label{eq:lambda-omega}
    \lambda_\omega = \frac{\varepsilon_\omega + 1}{2(\varepsilon_\omega - 1)}.
\end{align}
Moreover, $u_{2\omega}$ can be represented using the double-layer potentials $\mathcal{D}_{D}, \mathcal{D}_{\Omega_\epsilon}$ and the single-layer potentials $\mathcal{S}_{D}, \mathcal{S}_{\Omega_\epsilon}$ as follows:
\begin{align*}
    u_{2\omega}(x) = \mathcal{D}_{D} [\phi_i](x) + \mathcal{D}_{\Omega_\epsilon} [\phi_e](x) + \mathcal{S}_{D} [\psi_i](x) + \mathcal{S}_{ \Omega_\epsilon}[\psi_e](x), \quad x \in \mathbb{R}^2,
\end{align*}
where the quadruple $(\phi_i, \phi_e, \psi_i, \psi_e) \in L^2(\partial D) \times L^2(\partial \Omega_\epsilon) \times L_0^2(\partial D) \times L_0^2(\partial \Omega_\epsilon)$ satisfies
\begin{align*}
    \begin{cases}
    \ds \phi_i = 4\pi P_{\perp,i}^s & \textnormal{on } \partial D, \ms \\[8pt]
    \ds \phi_e = -4\pi P_{\perp,e}^s & \textnormal{on } \partial \Omega_\epsilon, \ms \\[8pt]
    \ds (\lambda_{2\omega} \mathcal{I} + \mathcal{K}_{ D}^*) [\psi_i] + \frac{\partial \mathcal{S}_{ \Omega_\epsilon} [\psi_e] }{\partial\nu_D}  = -\dfrac{\partial \mathcal{D}_{D} [\phi_i] }{ \partial \nu_{D}} - \dfrac{\partial \mathcal{D}_{ \Omega_\epsilon} [\phi_e] }{ \partial \nu_{D} } + \dfrac{ 4\pi \sigma_i^s}{ \varepsilon_{2\omega} - 1} & \textnormal{on } \partial D, \ms \\[8pt]
    \ds (\lambda_{2\omega} \mathcal{I} - \mathcal{K}_{ \Omega_\epsilon}^*) [\psi_e] - \frac{\partial \mathcal{S}_{D}[\psi_i] }{ \partial \nu_{ \Omega_\epsilon}}  =  \dfrac{ \partial \mathcal{ D}_{D} [\phi_i]}{\partial \nu_{\Omega_\epsilon}} + \dfrac{\partial \mathcal{D}_{\Omega_\epsilon} [\phi_e] }{\partial \nu_{\Omega_\epsilon} } -\dfrac{ 4\pi \sigma_e^s}{\varepsilon_{2\omega} - 1} & \textnormal{on } \partial \Omega_\epsilon.
    \end{cases}
\end{align*}
Here $\lambda_{2\omega}$ is given by
\begin{align}
    \label{eq:lambda-2omega}
    \lambda_{2\omega} = \frac{\varepsilon_{2\omega} + 1}{2(\varepsilon_{2\omega} - 1)}.
\end{align}
Furthermore, there exists a constant $C > 0$ such that
\begin{align*}
    \|\phi_i\|_{L^2(\partial D)} + \|\phi_e\|_{L^2(\partial \Omega_\epsilon)} + \|\psi_i\|_{L^2(\partial D)} + \|\psi_e\|_{L^2(\partial \Omega_\epsilon)} \leq C  \|\nabla H\|_{L^2(\partial D \cup \partial \Omega_\epsilon)} .
\end{align*}
\end{thm}
\begin{proof}
The proof is analogous to that presented in \cite{liu2026quantitative} and is therefore omitted.
\end{proof}
\section{Asymptotic analysis of the perturbed fields}
\label{sec:asymptotic analysis of the perturbed fields}
Let \(\nu_e\) and \(\tau_e\) denote the outward unit normal vector and the unit tangential vector on the reference boundary \(\partial \Omega\), respectively. For the perturbed outer boundary \(\partial \Omega_\epsilon\) defined in \eqref{eq:perturbed-boundary}, the corresponding outward unit normal \(\tilde{\nu}_{e}(\tilde{x})\) and unit tangential \(\tilde{\tau}_{e}(\tilde{x})\) on \(\partial \Omega_\epsilon\) admit the following asymptotic expansions as \(\epsilon \to 0\):
\begin{align*}
    %\label{eq:nu-expansion}
    \tilde{\nu}_{e} ( \tilde{x} ) & = \nu_{e}(x) - \epsilon \, f'(x) \, \tau_{e}(x) + O(\epsilon^2), \\[8pt]
    %\label{eq:tau-expansion}
    \tilde{\tau}_{e}(\tilde{x}) & = \tau_{e}(x) + \epsilon \, f'(x) \, \nu_{e}(x) + O(\epsilon^2),
\end{align*}
where \(f'(x)\) denotes the tangential derivative of \(f\) along \(\partial \Omega\), i.e., \(f' = \dfrac{\partial f}{\partial \tau_e}\).
\subsection{Derivation of asymptotic expansions via the field expansion method}
%\label{sec:asymptotic-expansions-via-field-expansion method}
With the above asymptotic relations for the boundary geometry, we are now in a position to derive the asymptotic expansions of the electric fields. By applying the field expansion method \cite{liu2026quantitative,nicholls2003analytic} and systematically matching the boundary conditions order by order in \(\epsilon\), we obtain the leading-order and first-order coupled systems. The results are summarized in the following theorem.

\begin{thm}
\label{thm-expansion}
Let $u_\omega$ and $u_{2\omega}$ be the solutions to \eqref{eq:governing-equation}. For sufficiently small $\epsilon >0 $ and for all $x \in \mathbb{R}^2$, the following pointwise asymptotic expansions hold
\begin{align}
    u_\omega(x) & = u_\omega^{(0)}(x) + \epsilon u_\omega^{(1)}(x) + O(\epsilon^2), \nonumber \\[8pt]
    \label{eq:u-2omega-expansion}
    u_{2\omega}(x) & = u_{2\omega}^{(0)}(x) + \epsilon u_{2\omega}^{(1)}(x) + O(\epsilon^2),
\end{align}
where the remainder $O(\epsilon^{2})$ depends only on the $\mathcal{C}^{2}$-norm of $\partial \Omega$ and $\mathcal{C}^{1}$-norm of $f$. $u_\omega^{(0)}$ and $u_{2\omega}^{(0)}$ are the solutions to the following leading-order coupled system,
\begin{align}
    \label{eq:leading-term-governing-equation}
    \renewcommand{\arraystretch}{1.5}
    \begin{array}{l l l}
    \left\{
    \begin{array}{l}
        \Delta u_\omega^{(0)} = 0 \\[8pt]
        u_\omega^{(0)}|_+ = u_\omega^{(0)}|_- \\[8pt]
        u_\omega^{(0)}|_+ = u_\omega^{(0)}|_- \\[8pt]
        \varepsilon_\omega \dfrac{\partial u_\omega^{(0)}}{\partial \nu_i} \bigg|_+ = \dfrac{\partial u_\omega^{(0)}}{\partial \nu_i} \bigg|_- \\[8pt]
        \dfrac{\partial u_\omega^{(0)}}{\partial \nu_e} \bigg|_+ = \varepsilon_\omega \dfrac{\partial u_\omega^{(0)}}{\partial \nu_e} \bigg|_- \\[8pt]
        u_\omega^{(0)} - H(x) = O(|x|^{-1})
    \end{array}
    \right.
    &
    \begin{array}{l}
    \Delta u_{2\omega}^{(0)} = 0 \\[8pt]
    u_{2\omega}^{(0)}|_+ - u_{2\omega}^{(0)}|_- = -4\pi P^{s,(0)}_{\perp,i} \\[8pt]
    u_{2\omega}^{(0)}|_+ - u_{2\omega}^{(0)}|_- = 4\pi P^{s,(0)}_{\perp,e} \\[8pt]
    \varepsilon_{2\omega} \dfrac{\partial u_{2\omega}^{(0)}}{\partial \nu_i} \bigg|_+ - \dfrac{\partial u_{2\omega}^{(0)}}{\partial \nu_i} \bigg|_- = 4\pi \sigma^{s,(0)}_{i} \\[8pt]
    \dfrac{\partial u_{2\omega}^{(0)}}{\partial \nu_e} \bigg|_+ - \varepsilon_{2\omega} \dfrac{\partial u_{2\omega}^{(0)}}{\partial \nu_e} \bigg|_- = -4\pi \sigma^{s,(0)}_{e} \\[8pt]
    u_{2\omega}^{(0)} = O(|x|^{-2})
    \end{array}
    &
    \begin{array}{l}
    \textnormal{in } \mathbb{R}^2 \setminus (\partial D \cup \partial \Omega), \\[8pt]
    \textnormal{on } \partial D, \\[8pt]
    \textnormal{on } \partial \Omega, \\[8pt]
    \textnormal{on } \partial D, \\[8pt]
    \textnormal{on } \partial \Omega, \\[8pt]
    \textnormal{as } |x| \to +\infty.
    \end{array}
    \end{array}
\end{align}
And the pair $(u^{(1)}_{\omega},u^{(1)}_{2\omega})$, is the unique solution to the following first-order coupled system
\begin{align}
    \label{eq:first-order-equation}
    \renewcommand{\arraystretch}{1.5}
    \begin{array}{l l l}
    \left\{
    \begin{array}{l}
        \Delta u^{(1)}_{\omega}=0 \\[8pt]
        u^{(1)}_{\omega}|_{+}-u^{(1)}_{\omega}|_{-}=0 \\[8pt]
        u^{(1)}_{\omega}|_{+}-u^{(1)}_{\omega}|_{-}=I_{1} \\[8pt]
        \varepsilon_{\omega} \dfrac{\partial u^{(1)}_{\omega}}{\partial\nu_i}\Big{|}_{+}- \dfrac{\partial u^{(1)}_{\omega}}{\partial\nu_i}\Big{|}_{-}=0 \\[8pt]
        \dfrac{\partial u^{(1)}_{\omega}}{\partial\nu_e}\Big{|}_{+}-\varepsilon_{\omega} \dfrac{\partial u^{(1)}_{\omega}}{\partial\nu_e}\Big{|}_{-}=I_{2} \\[8pt]
        u^{(1)}_{\omega}=O(|x|^{-1})
    \end{array}
    \right.
    &
    \begin{array}{l}
    \Delta u^{(1)}_{2\omega}=0 \\[8pt]
    u^{(1)}_{2\omega}|_{+}-u^{(1)}_{2\omega}|_{-}=- 4\pi P^{s,(1)}_{\perp,i} \\[8pt]
    u^{(1)}_{2\omega}|_{+}-u^{(1)}_{2\omega}|_{-}=I_{3} \\[8pt]
    \varepsilon_{2\omega} \dfrac{\partial u^{(1)}_{2\omega}}{\partial\nu_i}\Big{|}_{+}- \dfrac{\partial u^{(1)}_{2\omega}}{\partial\nu_i}\Big{|}_{-}=4\pi \sigma_i^{s,(1)}\\[8pt]
    \dfrac{\partial u^{(1)}_{2\omega}}{\partial\nu_e}\Big{|}_{+}-\varepsilon_{2\omega} \dfrac{\partial u^{(1)}_{2\omega}}{\partial\nu_e}\Big{|}_{-}=I_{4} \\[8pt]
    u^{(1)}_{2\omega}= o(1)
    \end{array}
    &
    \begin{array}{l}
    \textnormal{in } \mathbb{R}^2 \setminus (\partial D \cup \partial \Omega), \\[8pt]
    \textnormal{on } \partial D, \\[8pt]
    \textnormal{on } \partial \Omega, \\[8pt]
    \textnormal{on } \partial D, \\[8pt]
    \textnormal{on } \partial \Omega, \\[8pt]
    \textnormal{as } |x| \to +\infty,
    \end{array}
    \end{array}
\end{align}
with
\begin{align}
    \label{eq:boundary-term}
    \begin{cases}
    \ds I_{1} = f \left( \dfrac{\partial u^{(0)}_{\omega} }{\partial\nu_e}\Big{|}_{-} - \dfrac{\partial u^{(0)}_{\omega} }{ \partial\nu_e} \Big{|}_{+} \right) & \textnormal{on } \partial \Omega, \ms \\[8pt]
    \ds I_{2} = f \left( \varepsilon_{\omega} \dfrac{\partial^{2} u^{(0)}_{\omega}}{ \partial \nu_e^{2}} \Big{|}_{-} - \dfrac{ \partial^{2} u^{(0)}_{\omega} }{\partial\nu_e^{2}} \Big{|}_{+} \right) + f^{\prime} \left( \dfrac{ \partial u^{(0)}_{\omega} }{\partial \tau_e} \Big{|}_{+} - \varepsilon_{ \omega} \dfrac{\partial u^{(0)}_{\omega} }{\partial \tau_e} \Big{|}_{-} \right) & \textnormal{on } \partial \Omega, \ms \\[8pt]
    \ds I_{3} = f \left( \dfrac{\partial u^{(0)}_{ 2\omega} }{\partial\nu_e} \Big{|}_{-} - \dfrac{ \partial u^{(0)}_{2\omega} }{ \partial\nu_e} \Big{|}_{+}\right) + 4\pi P^{s,(1)}_{\perp,e} & \textnormal{on } \partial \Omega, \ms \\[8pt]
    \ds I_{4} = f \left( \varepsilon_{2\omega} \dfrac{ \partial^{2} u^{(0)}_{2\omega} }{\partial\nu_e^{2} } \Big{|}_{-} - \dfrac{\partial^{2} u^{(0)}_{ 2\omega} }{\partial\nu_e^{2}} \Big{|}_{+} \right) + f^{\prime} \left( \dfrac{\partial u^{(0)}_{2\omega} }{ \partial \tau_e} \Big{|}_{+} - \varepsilon_{ 2\omega} \dfrac{\partial u^{(0)}_{2\omega} }{ \partial \tau_e} \Big{|}_{-}\right) - 4\pi \sigma^{s,(1)}_{e} & \textnormal{on } \partial \Omega,
    \end{cases}
\end{align}
where the quantities $P^{s,(0)}_{\perp,e}$, $P^{s,(1)}_{\perp,e}$, $P^{s,(0)}_{\perp,i}$, $P^{s,(1)}_{\perp,i}$ are defined in \eqref{eq:Ps-perp-e-leading}, \eqref{eq:Ps-perp-e-first}, \eqref{eq:Ps-perp-i-leading}, \eqref{eq:Ps-perp-i-first}, respectively; and the surface charge terms $\sigma_e^{s,(0)}$, $\sigma_e^{s,(1)}$ and $\sigma_i^{s,(0)}$, $\sigma_i^{s,(1)}$ are given by \eqref{eq:sigma-s-e-expansion} and \eqref{eq:sigma-s-i-expansion}, respectively.
\end{thm}
\begin{proof}
We assume the pointwise expansion for \( u_{\omega} \)
\begin{align*}
    u_\omega(x)  = u_\omega^{(0)}(x) + \epsilon u_\omega^{(1)}(x) + O(\epsilon^2).
\end{align*}

For $x\in\partial \Omega$, the perturbed boundary point is $\tilde{x}=x + \epsilon f(x)\nu_e(x) \in \partial\Omega_\epsilon$. Then we have the following Taylor expansions
\begin{align}
    \label{eq:u-omega-Omega-e}
    u_\omega |_\pm(\tilde{x}) & = u_\omega^{(0)}|_\pm(x) + \epsilon \left(u_\omega^{(1)} + f \frac{\partial u_\omega^{(0)} }{\partial\nu_e} \right) \Big|_\pm(x) + O (\epsilon^2), \\[4pt]
    \label{eq:u-omega-normal-derivative-Omega-e}
    \frac{\partial u_\omega}{\partial \tilde{\nu}_e} \Big|_\pm(\tilde{x}) & = \frac{\partial u_\omega^{(0)}}{\partial\nu_e} + \epsilon\left(\frac{\partial u_\omega^{(1)}}{\partial\nu_e} + f \frac{\partial^2 u_\omega^{(0)}}{\partial\nu_e^2} - f' \frac{\partial u_\omega^{(0)} }{\partial\tau_e} \right) \Big|_\pm(x) + O (\epsilon^2), \\[4pt]
    \label{eq:u-omega-tangential-derivative-Omega-e}
    \frac{\partial u_\omega}{\partial\tilde{\tau}_e} \Big|_{-} (\tilde{x}) & = \frac{\partial u_\omega^{(0)} }{\partial \tau_e} \Big|_{-} (x) + \epsilon \left( \frac{\partial u_\omega^{(1)} }{\partial \tau_e} + f \frac{\partial}{\partial \nu_e} \left( \frac{\partial u_\omega^{(0)}} {\partial \tau_e} \right) + f' \frac{\partial u_\omega^{(0)} }{\partial \nu_e} \right)\Big|_{-} (x) + O \left( \epsilon^2 \right).
\end{align}
By using the transmission conditions on $\partial \Omega_\epsilon$, we deduce from \eqref{eq:u-omega-Omega-e} and \eqref{eq:u-omega-normal-derivative-Omega-e} that
\begin{align*}
    u_\omega^{(0)}|_{+} = u_\omega^{(0)}|_{-} \quad & \text{on }  \partial \Omega, \\[8pt]
    \frac{\partial u^{(0)}_\omega}{\partial \nu_e} \Big|_{+} = \varepsilon_\omega \frac{\partial u^{(0)}_\omega}{\partial \nu_e} \Big|_{-} \quad & \text{on }  \partial \Omega,
\end{align*}
and
\begin{align*}
    u_\omega^{(1)}|_{+} - u_\omega^{(1)}|_{-} = f \left( \frac{ \partial u_\omega^{(0)} }{ \partial \nu_e} \Big|_{-} - \frac{\partial u_\omega^{(0)} }{ \partial \nu_e} \Big|_{+} \right) \quad & \text{on }  \partial \Omega, \\[4pt]
    \frac{\partial u_\omega^{(1)} }{\partial \nu_e} \Big|_{+} - \varepsilon_\omega \frac{ \partial u_\omega^{(1)} }{\partial \nu_e} \Big|_{-} = f \left( \varepsilon_\omega \frac{\partial^2 u_\omega^{(0)} }{\partial \nu_e^2} \Big|_{-} - \frac{\partial^2 u_\omega^{(0)} }{\partial \nu_e^2} \Big|_{+} \right) + f' \left(\frac{\partial u_\omega^{(0)} }{\partial \tau_e} \Big|_{+} - \varepsilon_\omega \frac{\partial u_\omega^{(0)} }{\partial \tau_e} \Big|_{-} \right) \quad & \text{on }  \partial \Omega.
\end{align*}

For $x\in\partial D$, from the field expansions, we obtain
\begin{align}
    \label{eq:u-omega-Di}
    u_\omega|_\pm(x) & = u_\omega^{(0)}|_\pm(x) + \epsilon u_\omega^{(1)}|_\pm(x) + O(\epsilon^2), \\[8pt]
    \label{eq:u-omega-normal-derivative-Di}
    \frac{\partial  u_\omega}{\partial\nu_i} \Big{|}_\pm(x) & = \frac{\partial u_\omega^{(0)} }{\partial\nu_i} \Big{|}_\pm(x)+ \epsilon \frac{\partial u_\omega^{(1)} }{\partial\nu_i} \Big{|}_\pm(x)+O(\epsilon^2), \\[4pt]
    \label{eq:u-omega-tangential-derivative-Di}
    \frac{\partial u_\omega}{\partial\tau_i} \Big{|}_+(x) & = \frac{\partial u_\omega^{(0)} }{\partial\tau_i} \Big{|}_+(x) + \epsilon \frac{\partial u_\omega^{(1)} }{\partial\tau_i} \Big{|}_+(x) + O(\epsilon^2).
\end{align}
By using the transmission conditions on $\partial D$, we deduce from \eqref{eq:u-omega-Di} and \eqref{eq:u-omega-normal-derivative-Di} that
\begin{align*}
    u_\omega^{(0)}|_{+} = u_\omega^{(0)}|_{-} \quad & \text{on }  \partial D, \\[8pt]
    \varepsilon_\omega \frac{\partial u^{(0)}_\omega}{\partial \nu_i} \Big|_{+} = \frac{\partial u^{(0)}_\omega}{\partial \nu_i} \Big|_{-} \quad & \text{on }  \partial D,
\end{align*}
and
\begin{align*}
    u_\omega^{(1)}|_{+} = u_\omega^{(1)}|_{-} \quad & \text{on }  \partial D, \\[8pt]
    \varepsilon_\omega \frac{\partial u_\omega^{(1)}}{\partial \nu_i} \Big|_{+} = \frac{\partial u_\omega^{(1)}}{\partial \nu_i} \Big|_{-} \quad & \text{on }  \partial D.
\end{align*}

To perform a Taylor expansion of the surface charge functions \(\sigma_e^s\) and \(\sigma_i^s\) defined by \eqref{eq:sigma-s-e} and \eqref{eq:sigma-s-i}, respectively, we parameterize the perturbed outer boundary \(\partial\Omega_\epsilon\) as \(X_e(t)=(x_1^{(e)}(t),x_2^{(e)}(t))\) and the unperturbed inner boundary \(\partial D\) as \(X_i(t)=(x_1^{(i)}(t),x_2^{(i)}(t))\). For the outer boundary, we have
\begin{align}
    \label{eq:sigma-s-e-parameterize}
    \sigma_e^s = -\frac{1}{h_e} \frac{dP_{\parallel,e}^s}{dt},
\end{align}
where \(h_e\) is the Lam\'e coefficient of \(\partial\Omega_\epsilon\). Using perturbation theory, the Lam\'e coefficient expands as
\begin{align}
    \label{eq:Lame-expansion}
    h_e = h_e^{(0)} + \epsilon\,f h_e^{(0)}\kappa_e + O(\epsilon^2),
\end{align}
with \(\kappa_e\) the curvature of the reference boundary \(\partial \Omega\). For the inner boundary, which remains unperturbed, we simply have
\begin{align}
    \label{eq:sigma-s-i-parameterize}
    \sigma_i^s = -\frac{1}{h_i} \frac{ dP_{\parallel,i}^s}{dt}, \qquad h_i = r_i,
\end{align}
where \(h_i\) is the Lam\'e coefficient of \(\partial D\), and no further expansion is needed.

For $x\in\partial \Omega$, let $\tilde{x}=x+ \epsilon f(x)\nu_e(x) \in \partial \Omega_\epsilon$. From \eqref{eq:Pse-perp}, \eqref{eq:Pse-parallel}, \eqref{eq:u-omega-normal-derivative-Omega-e} and \eqref{eq:u-omega-tangential-derivative-Omega-e}, it follows that
\begin{align}
    \label{eq:Ps-perp-e-expansion}
    P^{s}_{\perp,e}(\tilde{x}) & = P^{s,(0)}_{\perp,e} + \epsilon P^{s,(1)}_{\perp,e} + O\left(\epsilon^2\right),\\[4pt]
    \label{eq:Ps-parallel-e-expansion}
    P^{s}_{\parallel,e}(\tilde{x}) & =P^{s,(0)}_{\parallel,e} + \epsilon P^{s,(1)}_{\parallel,e} + O\left(\epsilon^2\right),
\end{align}
where
\begin{align}
    \label{eq:Ps-perp-e-leading}
    P^{s,(0)}_{\perp,e} & = \chi_{\perp}^s \varepsilon_\omega^2 \left( \frac{\partial u^{(0)}_\omega }{\partial \nu_e} \Big|_{-}(x) \right)^2, \\[4pt]
    \label{eq:Ps-perp-e-first}
    P^{s,(1)}_{\perp,e} & = 2 \chi_{\perp}^s \varepsilon_\omega^2  \left( \frac{\partial u^{(0)}_\omega}{\partial \nu_e} \Big|_{-}(x) \right) \left( \frac{\partial u_\omega^{(1)} }{ \partial \nu_e} \Big|_{-}(x) + f \frac{\partial^2 u_\omega^{(0)} }{ \partial \nu_e^2} \Big|_{-}(x) - f' \frac{ \partial u_\omega^{(0)} }{ \partial \tau_e } \Big|_{-}(x) \right),\\[4pt]
    \label{eq:Ps-parallel-e-leading}
    P^{s,(0)}_{\parallel,e} & = 2 \chi_{\parallel}^s \varepsilon_\omega  \left( \frac{\partial u^{(0)}_\omega}{ \partial \nu_e} \Big|_{-}(x)\right)  \left( \frac{\partial u^{(0)}_\omega}{\partial \tau_e} \Big|_{-}(x) \right), \\[4pt]
    P^{s,(1)}_{\parallel,e} & = 2 \chi_{\parallel}^s \varepsilon_\omega \left( \frac{\partial u^{(0)}_\omega}{ \partial \nu_e} \Big|_{-}(x)\right) \left( \frac{ \partial u_\omega^{(1)} }{ \partial \tau_e} \Big|_{-}(x) + f \frac{\partial}{\partial \nu_e} \left( \frac{\partial u_\omega^{(0)}}{ \partial \tau_e}\right) \Big|_{-}(x) + f' \frac{\partial u_\omega^{(0)} }{\partial \nu_e}\Big|_{-}(x) \right) \nonumber\\[4pt]
    \label{eq:Ps-parallel-e-first}
    & \quad + 2 \chi_{\parallel}^s \varepsilon_\omega \left( \frac{\partial u^{(0)}_\omega }{ \partial \tau_e} \Big|_{-}(x) \right) \left( \frac{\partial u_\omega^{(1)} }{ \partial \nu_e} \Big|_{-}(x) + f \frac{\partial^2 u_\omega^{(0)} }{ \partial \nu_e^2} \Big|_{-}(x) - f' \frac{\partial u_\omega^{(0)} }{\partial \tau_e}\Big|_{-}(x) \right).
\end{align}
By \eqref{eq:sigma-s-e-parameterize}, \eqref{eq:Lame-expansion} and \eqref{eq:Ps-parallel-e-expansion}, we have
\begin{align}
    \sigma_e^s & = -\frac{1}{h_e^{(0)}} \frac{\d P^{s,(0)}_{\parallel,e}}{\d t} + \epsilon \left( \frac{f\kappa_e}{h_e^{(0)}} \frac{\d P^{s,(0)}_{ \parallel,e}}{\d t} - \frac{1}{h_e^{(0)}} \frac{\d P^{s,(1)}_{\parallel,e}}{\d t}  \right) + O \left( \epsilon^2 \right) \nonumber \\[8pt]
    \label{eq:sigma-s-e-expansion}
    &= \sigma_e^{s,(0)} + \epsilon \sigma_e^{s,(1)} + O \left( \epsilon^2 \right).
\end{align}

For $x\in\partial D$, from \eqref{eq:Psi-perp}, \eqref{eq:Psi-parallel}, \eqref{eq:u-omega-normal-derivative-Di} and \eqref{eq:u-omega-tangential-derivative-Di}, it follows that
\begin{align}
    \label{eq:Ps-perp-i-expansion}
    P^{s}_{ \perp,i }(x) & = P^{s,(0)}_{\perp,i} + \epsilon P^{s,(1)}_{\perp,i} + O \left( \epsilon^2 \right),\\[8pt]
    \label{eq:Ps-parallel-i-expansion}
    P^{s}_{ \parallel,i }(x) & =P^{s,(0)}_{ \parallel,i } + \epsilon P^{s,(1)}_{ \parallel,i } + O \left( \epsilon^2 \right),
\end{align}
where
\begin{align}
    \label{eq:Ps-perp-i-leading}
    P^{s,(0)}_{\perp,i} & = \chi_{\perp}^s \varepsilon_{\omega}^2 \left( \frac{\partial u_{\omega}^{(0)} }{ \partial \nu_i} \Big|_{+}(x) \right)^2 , \\[4pt]
    \label{eq:Ps-perp-i-first}
    P^{s,(1)}_{\perp,i} & = 2\chi_{\perp}^s \varepsilon_{\omega}^2  \left( \frac{\partial u_{\omega}^{(0)} }{ \partial \nu_i } \Big|_{+}(x) \right) \left( \frac{\partial u_{\omega}^{(1)} }{ \partial \nu_i } \Big|_{+}(x) \right) , \\[4pt]
    \label{eq:Ps-parallel-i-leading}
    P^{s,(0)}_{\parallel,i} & = 2 \chi_{\parallel}^s \varepsilon_{\omega} \left( \frac{ \partial u_{\omega}^{(0)} }{ \partial \nu_i } \Big|_{+}(x) \right) \left( \frac{ \partial u_{\omega}^{(0)} }{ \partial \tau_i} \Big|_{+}(x) \right ) , \\[4pt]
    \label{eq:Ps-parallel-i-first}
    P^{s,(1)}_{\parallel,i} & = 2\chi_{\parallel}^s \varepsilon_{\omega}  \left[ \left( \frac{\partial u_{\omega}^{(0)} }{ \partial \nu_i} \Big|_{+}(x) \right) \left( \frac{ \partial u_{\omega}^{(1)}}{\partial \tau_i} \Big|_{+}(x) \right) +  \left( \frac{ \partial u_{\omega}^{(0)} }{ \partial \tau_i } \Big|_{+}(x) \right) \left( \frac{ \partial u_{\omega}^{(1)} }{ \partial \nu_i } \Big|_{+}(x) \right) \right].
\end{align}
By \eqref{eq:sigma-s-i-parameterize} and \eqref{eq:Ps-parallel-i-expansion}, we have
\begin{align}
    \sigma_i^s & = -\frac{1}{h_i} \frac{\d P^{s,(0)}_{\parallel,i}}{\d t} - \epsilon \frac{1}{h_i} \frac{\d P^{s,(1)}_{\parallel,i}}{\d t} + O\left(\epsilon^2\right) \nonumber \\[8pt]
    \label{eq:sigma-s-i-expansion}
    & = \sigma_i^{s,(0)} + \epsilon \sigma_i^{s,(1)} + O \left( \epsilon^2 \right).
\end{align}

Similarly, we assume the following pointwise expansion for \( u_{2\omega} \)
\begin{align*}
    u_{2\omega}(x)  = u_{2\omega}^{(0)}(x) + \epsilon u_{2\omega}^{(1)}(x) + O(\epsilon^2).
\end{align*}
For $x\in\partial \Omega$, let $\tilde{x}=x+ \epsilon f(x)\nu_e(x) \in \partial\Omega_\epsilon$. Following the same procedure as in \eqref{eq:u-omega-Omega-e} and \eqref{eq:u-omega-normal-derivative-Omega-e}, we obtain the following Taylor expansions
\begin{align}
    \label{eq:u-2omega-Omega-e}
    u_{2\omega}|_\pm (\tilde{x}) & = u_{2\omega}^{(0) }|_\pm(x) + \epsilon \left( u_{2\omega}^{(1)} + f \frac{ \partial u_{2\omega}^{(0)} }{ \partial\nu_e} \right) \Big|_\pm(x) + O(\epsilon^2), \\[4pt]
    \label{eq:u-2omega-normal-derivative-Omega-e}
    \frac{ \partial u_{2\omega} }{ \partial \tilde{\nu}_e } \Big|_\pm (\tilde{x}) & = \frac{\partial u_{2\omega}^{(0)} }{ \partial\nu_e } \Big|_\pm(x) + \epsilon \left( \frac{ \partial u_{2\omega}^{(1)} }{ \partial\nu_e} + f \frac {\partial^2 u_{2\omega}^{(0)} }{ \partial \nu_e^2} - f' \frac{ \partial u_{2\omega}^{(0)} }{ \partial \tau_e} \right) \Big|_\pm(x) + O(\epsilon^2).
\end{align}
From \eqref{eq:Ps-perp-i-expansion}, \eqref{eq:sigma-s-e-expansion}, \eqref{eq:u-2omega-Omega-e} and \eqref{eq:u-2omega-normal-derivative-Omega-e}, the transmission conditions on $\partial\Omega_\epsilon$ immediately yield
\begin{align*}
    u^{(0)}_{2\omega}|_{+} - u^{(0)}_{2\omega}|_{-} = 4\pi P^{s,(0)}_{\perp,e}  & \quad \text{on } \partial \Omega, \\[8pt]
    \frac{ \partial u^{(0)}_{2\omega} }{ \partial \nu_e} \Big|_{+} - \varepsilon_{ 2\omega} \frac{\partial u^{(0)}_{ 2\omega} }{\partial \nu_e} \Big|_{-} = -4\pi \sigma_e^{s,(0)}  & \quad \text{on } \partial \Omega,
\end{align*}
and
\begin{align*}
    u_{2\omega}^{(1)}|_{+} - u_{2\omega}^{(1)}|_{-} = f \left( \frac{ \partial u_{2\omega}^{(0)} }{ \partial \nu_e} \Big|_{-} - \frac{\partial u_{2\omega}^{(0)} }{ \partial \nu_e} \Big|_{+}\right) + 4\pi P^{s,(1)}_{\perp,e} &\quad \text{on } \partial \Omega, \\[4pt]
    \frac{ \partial u_{2\omega}^{(1)} }{\partial \nu_e} \Big|_{+} - \varepsilon_{2\omega} \frac{\partial u_{2\omega}^{(1)} }{ \partial \nu_e} \Big|_{-} = f \left( \varepsilon_{2\omega} \frac{ \partial^2 u_{2\omega}^{(0)} }{ \partial \nu_e^2 } \Big|_{-} - \frac{ \partial^2 u_{ 2\omega }^{(0)} }{ \partial \nu_e^2} \Big|_{+} \right) + f'\left(\frac{\partial u_{2\omega}^{(0)} }{ \partial \tau_e} \Big|_{+} - \varepsilon_{2\omega} \frac{ \partial u_{2\omega }^{(0)}}{\partial \tau_e} \Big|_{-} \right) - 4\pi \sigma_e^{s,(1)}  & \quad \text{on } \partial \Omega.
\end{align*}

For $x\in\partial D$, by applying the field expansions and following the same procedure as in \eqref{eq:u-omega-Di} and \eqref{eq:u-omega-normal-derivative-Di}, we obtain
\begin{align*}
    %\label{eq:u-2omega-Di}
    u_{2\omega}|_\pm(x) & = u_{2\omega}^{(0)}|_\pm(x) + \epsilon u_{2\omega}^{(1)}|_\pm(x) + O(\epsilon^2), \\[8pt]
    \label{eq:u-2omega-normal-derivative-Di}
    \frac{\partial u_{2\omega} }{ \partial\nu_i } \Big{|}_\pm(x) & = \frac{\partial u_{2\omega}^{(0)} }{ \partial\nu_i} \Big{|}_\pm(x) + \epsilon \frac{\partial  u_{2\omega}^{(1)} }{ \partial \nu_i } \Big{|}_\pm(x) + O(\epsilon^2).
\end{align*}
From \eqref{eq:Ps-perp-e-expansion}, \eqref{eq:sigma-s-i-expansion}, \eqref{eq:u-2omega-Omega-e} and \eqref{eq:u-2omega-normal-derivative-Omega-e}, the transmission conditions on $\partial D$ immediately yield
\begin{align*}
    u^{(0)}_{2\omega}|_{+} - u^{(0)}_{2\omega}|_{-} = -4\pi P^{s,(0)}_{\perp,i}  & \quad \text{on } \partial D, \\[8pt]
    \varepsilon_{2\omega} \frac{\partial u^{(0)}_{ 2\omega} }{ \partial \nu_i} \Big|_{+} - \frac{ \partial u^{(0)}_{2\omega} }{ \partial \nu_i} \Big|_{-} = 4 \pi \sigma_i^{s,(0)}  & \quad \text{on } \partial D,
\end{align*}
and
\begin{align*}
    u_{2\omega}^{(1)}|_{+} - u_{2\omega}^{(1)}|_{-} = - 4\pi P^{s,(1)}_{\perp,i} & \quad \text{on } \partial D, \\[8pt]
    \varepsilon_{2\omega} \frac{ \partial u_{ 2\omega }^{(1)} }{ \partial \nu_i } \Big|_{+} - \frac{ \partial u_{2\omega}^{(1)} }{ \partial \nu_i} \Big|_{-} = 4\pi \sigma_i^{s,(1)}  & \quad \text{on } \partial D.
\end{align*}
\end{proof}

\subsection{Layer-potential representations of the asymptotic solutions}
\label{sec:layer-potential-asymptotic-solutions}
In this subsection, we establish the representation formulas for the solutions to the leading-order and first-order coupled systems in Theorem \ref{thm-expansion}. Using the layer-potential theory introduced in Section \ref{sec:layer-potential-formulation-and-well-posedness}, we express the leading-order and first-order solutions in terms of single-layer and double-layer potentials.

Applying the same layer-potential argument as in Theorem~\ref{well-posedness}, the fundamental-field component \(u_\omega^{(0)}\) of the leading-order system \eqref{eq:leading-term-governing-equation} admits the representation
\vspace{6pt}
\begin{equation}
    \label{eq:u-omega-leading}
    u_\omega^{(0)}(x) = H(x) + \mathcal{S}_{D}[\varphi_i^{(0)}](x) + \mathcal{S}_{\Omega}[\varphi_e^{(0)}](x), \quad x \in \mathbb{R}^2.
\vspace{6pt}
\end{equation}
Here the density pair \( \bigl( \varphi_i^{(0)}, \varphi_e^{(0)} \bigr) \in L_0^2( \partial D) \times L_0^2(\partial\Omega) \) is determined by the coupled boundary integral system
\begin{align}
    \label{eq:u-omega-leading-density}
    \begin{cases}
    \ds ( \lambda_\omega \mathcal{I} + \mathcal{K}_{D}^* ) [\varphi_i^{(0)}] + \frac{ \partial \mathcal{S}_{\Omega} [\varphi_e^{(0)}] }{ \partial\nu_i}  = -\dfrac{ \partial H }{ \partial \nu_i} & \text{on } \partial D, \\[8pt]
    \ds ( \lambda_\omega \mathcal{I} - \mathcal{K}_{\Omega}^* ) [\varphi_e^{(0)}] - \frac{ \partial \mathcal{S}_D  [\varphi_i^{(0)}] }{ \partial\nu_e } = \dfrac{ \partial H }{ \partial \nu_e } & \text{on } \partial \Omega.
    \end{cases}
\end{align}
The parameter \( \lambda_\omega \) is defined in \eqref{eq:lambda-omega}. Similarly, the second-harmonic component \(u_{2\omega}^{(0)}\) of \eqref{eq:leading-term-governing-equation} admits the representation
\begin{align}
    \label{eq:u-2omega-leading}
    u_{2\omega}^{(0)}(x) = \mathcal{D}_{D}[\phi_i^{(0)}](x) + \mathcal{D}_{\Omega}[\phi_e^{(0)}](x) + \mathcal{S}_{D}[\psi_i^{(0)}](x) + \mathcal{S}_{\Omega}[\psi_e^{(0)}](x), \quad x \in \mathbb{R}^2.
\end{align}
The corresponding density vector \(\bigl(\phi_i^{(0)}, \phi_e^{(0)}, \psi_i^{(0)}, \psi_e^{(0)}\bigr) \in L^2(\partial D) \times L^2(\partial \Omega) \times L_0^2(\partial D) \times L_0^2(\partial \Omega)\) is determined by
\begin{align}
    \label{eq:u-2omega-leading-density}
    \begin{cases}
    \ds \phi_i^{(0)} = 4\pi P_{\perp,i}^{s,(0)} & \text{on } \partial D, \\[8pt]
    \ds \phi_e^{(0)} = -4\pi P_{\perp,e}^{s,(0)} & \text{on } \partial \Omega, \\[8pt]
    \ds (\lambda_{2\omega} \mathcal{I} + \mathcal{K}_{D}^*) [\psi_i^{(0)}] +  \frac{\partial \mathcal{S}_{\Omega} [\psi_e^{(0)}] }{\partial\nu_i}   = - \dfrac{\partial \mathcal{D}_{D} [\phi_i^{(0)}] }{ \partial \nu_i} - \dfrac{\partial \mathcal{D}_{\Omega} [\phi_e^{(0)}] }{ \partial \nu_i} + \dfrac{ 4\pi \sigma_i^{s,(0)}} {\varepsilon_{2\omega} - 1} & \text{on } \partial D, \\[8pt]
    \ds ( \lambda_{2\omega} \mathcal{I} - \mathcal{K}_{\Omega}^*) [\psi_e^{(0)}] - \frac{ \partial \mathcal{S}_D [\psi_i^{(0)}] }{ \partial\nu_e}  =  \dfrac{\partial \mathcal{D}_{D} [\phi_i^{(0)}] }{ \partial \nu_e} + \dfrac{ \partial \mathcal{D}_{\Omega} [\phi_e^{(0)}] }{ \partial \nu_e} - \dfrac{4\pi \sigma_e^{s,(0)} }{ \varepsilon_{2\omega} - 1} & \text{on } \partial \Omega,
    \end{cases}
\end{align}
where the parameter \(\lambda_{2\omega}\) is defined in
\eqref{eq:lambda-2omega}.

For the first-order system \eqref{eq:first-order-equation}, the fundamental-field correction \(u_\omega^{(1)}\) admits the representation
\begin{align}
    \label{eq:u-omega-first}
    u_\omega^{(1)}(x) = \mathcal{D}_{\Omega} [\vartheta_e^{(1)}](x) + \mathcal{S}_{D} [\varphi_i^{(1)}](x) + \mathcal{S}_{\Omega} [\varphi_e^{(1)}](x), \quad x \in \mathbb{R}^2.
\end{align}
The density triplet \(\bigl(\vartheta_e^{(1)}, \varphi_i^{(1)}, \varphi_e^{(1)}\bigr) \in L_0^2(\partial \Omega) \times L_0^2(\partial D) \times L_0^2(\partial \Omega)\) is determined by
\begin{align}
    \label{eq:u-omega-first-density}
    \begin{cases}
    \ds \vartheta_e^{(1)} = -I_1 & \text{on } \partial \Omega, \\[8pt]
    \ds ( \lambda_\omega \mathcal{I} + \mathcal{K}_{D}^*) [\varphi_i^{(1)}] + \frac{ \partial \mathcal{S}_{\Omega} [\varphi_e^{(1)}] }{\partial\nu_i} = -\dfrac{\partial \mathcal{D}_{\Omega} [\vartheta_e^{(1)}]}{\partial \nu_i} & \text{on } \partial D, \\[8pt]
    \ds (\lambda_\omega \mathcal{I} - \mathcal{K}_{\Omega}^*) [\varphi_e^{(1)}] - \frac{ \partial \mathcal{S}_D [\varphi_i^{(1)}] }{ \partial\nu_e} = \dfrac{\partial \mathcal{D}_{\Omega} [\vartheta_e^{(1)}] }{ \partial \nu_e} + \dfrac{I_2} {\varepsilon_\omega - 1} & \text{on } \partial \Omega.
    \end{cases}
\end{align}

The second-harmonic correction \(u_{2\omega}^{(1)}\) admits the representation
\begin{align}
    \label{eq:u-2omega-first}
    u_{2\omega}^{(1)}(x) = \mathcal{D}_{D} [\phi_i^{(1)}](x) + \mathcal{D}_{\Omega} [\phi_e^{(1)}](x) + \mathcal{S}_{D} [\psi_i^{(1)}](x) + \mathcal{S}_{\Omega} [\psi_e^{(1)}](x), \quad x \in \mathbb{R}^2.
\end{align}
where the corresponding density vector  \(\bigl(\phi_i^{(1)}, \phi_e^{(1)}, \psi_i^{(1)}, \psi_e^{(1)}\bigr) \in L_0^2(\partial D) \times L_0^2(\partial \Omega) \times L_0^2(\partial D) \times L_0^2(\partial \Omega)\) satisfies
\begin{align}
    \label{eq:u-2omega-first-density}
    \begin{cases}
    \ds \phi_i^{(1)} = 4\pi P^{s,(1)}_{\perp,i} & \text{on } \partial D, \\[8pt]
    \ds \phi_e^{(1)} = -I_{3} & \text{on } \partial \Omega, \\[8pt]
    \ds (\lambda_{2\omega} \mathcal{I} + \mathcal{K}_{D}^*) [\psi_i^{(1)}] + \frac{\partial \mathcal{S}_{\Omega} [\psi_e^{(1)}]} {\partial\nu_i}  = - \dfrac{\partial \mathcal{D}_{D} [\phi_i^{(1)}] }{\partial \nu_i} - \dfrac{\partial \mathcal{D}_{\Omega} [\phi_e^{(1)}] }{\partial \nu_i} + \dfrac{ 4\pi \sigma_i^{s,(1)} }{ \varepsilon_{2\omega} - 1} & \text{on } \partial D, \\[8pt]
    \ds (\lambda_{2\omega} \mathcal{I} - \mathcal{K}_{\Omega}^*) [\psi_e^{(1)}] - \frac{\partial \mathcal{S}_D [\psi_i^{(1)}]} {\partial\nu_e}  = \dfrac{\partial \mathcal{D}_{D} [\phi_i^{(1)}]} {\partial \nu_e} + \dfrac{\partial \mathcal{D}_\Omega [\phi_e^{(1)}] } { \partial \nu_e } + \dfrac{ I_{4} } { \varepsilon_{2\omega} - 1} & \text{on } \partial \Omega.
\end{cases}
\end{align}

\section{Proofs of the main results}
\label{sec:proofs of the main results}
In this section, we consider the second harmonic generation in a hollow nanowire under the uniform background field specified above. The cross-section of the hollow nanowire consists of two concentric boundaries: the inner boundary \(\partial D\) is a circle of radius \(r_i\) centered at the origin, and the outer boundary \(\partial \Omega_\epsilon\) is a perturbation of the reference circle \(\partial\Omega\) of radius \(r_e\). Throughout this section, we set \( D = \{x\in\mathbb{R}^2:|x|<r_i\} \) and \( \Omega = \{x\in\mathbb{R}^2:|x|<r_e\} \). Then, from \eqref{eq:perturbed-boundary}, the perturbed outer boundary \(\partial \Omega_\epsilon\) in polar coordinates is expressed as
\begin{align*}
    \begin{cases}
    \tilde{x}_{e}^{(1)} = (r_e + \epsilon f(\theta)) \cos \theta, \\[6pt]
    \tilde{x}_{e}^{(2)} = (r_e + \epsilon f(\theta)) \sin \theta.
    \end{cases}
\end{align*}
For the shape function considered in Theorem~\ref{thm:mainthm2}, we have
\begin{align}
    \label{eq:f-series}
    f(\theta) = r_e \cos(n\theta).
\end{align}
The inner boundary remains unperturbed and is described by
\begin{align*}
    x_{i}^{(1)} = r_i \cos \theta, \qquad x_{i}^{(2)} = r_i \sin \theta.
\end{align*}

We begin by proving Theorem~\ref{thm:mainthm1} for the unperturbed configuration.
\begin{proof}[Proof of Theorem \ref{thm:mainthm1}]
From the layer-potential theory in Subsection \ref{sec:layer-potential-asymptotic-solutions}, the leading-order fundamental field $u_\omega^{(0)}$ can be expressed as \eqref{eq:u-omega-leading}. For the concentric disks $D$ and $\Omega$, using \eqref{eq:single-layer-explicit}, \eqref{eq:single-layer-derivative-explicit} and \eqref{eq:NP-operator-explicit} and solving the coupled boundary integral equations \eqref{eq:u-omega-leading-density}, we obtain the density functions
\begin{align*}
    %\label{eq:u-omega-leading-density-explicit}
    \begin{cases}
    \ds \varphi_i^{(0)} = -E \frac{-4 \varepsilon_{ \omega-1}}{\varepsilon_{\omega+1}^2 - \gamma^2 \varepsilon_{\omega-1}^2} \cos \theta & \text{on } \partial D, \ms \\[4pt]
    \ds \varphi_e^{(0)} = -E \frac{2\varepsilon_{ \omega-1} (\varepsilon_{\omega+1} - \gamma^2 \varepsilon_{\omega-1})}{\varepsilon_{\omega+1}^2 - \gamma^2 \varepsilon_{\omega-1}^2} \cos \theta & \text{on } \partial \Omega.
    \end{cases}
\end{align*}
Substituting the densities $(\varphi_i^{(0)}, \varphi_e^{(0)})$ into the representation formula \eqref{eq:u-omega-leading} and using \eqref{eq:single-layer-explicit}, we obtain the following explicit expression for \(u_\omega^{(0)}\)
\begin{align}
    \label{eq:u-omega-leading-explicit}
    u_\omega^{(0)} =
    \begin{cases}
    \ds  -E \alpha_{i}^{(0)} r \cos \theta, &\quad r < r_i, \\[8pt]
    \ds -E \alpha_{m-}^{(0)} r\cos \theta -E \alpha_{m+}^{(0)} \dfrac{1}{r}\cos \theta, &\quad r_i < r < r_e, \\[8pt]
    \ds  -Er \cos \theta -E \alpha_{e}^{(0)} \dfrac{1}{r} \cos \theta, & \quad r > r_e.
    \end{cases}
\end{align}
where
\begin{align*}
    \ds \alpha_{i}^{(0)} & = \dfrac{ 4\varepsilon_\omega}{ \varepsilon_{\omega+1}^2 - \gamma^2 \varepsilon_{\omega-1}^2}, \qquad  \alpha_{e}^{(0)}  =  \dfrac{ \varepsilon_{ \omega+1} \varepsilon_{\omega-1} }{ \varepsilon_{\omega+1}^2 - \gamma^2 \varepsilon_{\omega-1}^2} (r_i^2-r_e^2), \\[8pt]
    \alpha_{m-}^{(0)} & = \dfrac{ 2\varepsilon_{ \omega+1} }{ \varepsilon_{\omega+1}^2 - \gamma^2 \varepsilon_{\omega-1}^2 }, \qquad \alpha_{m+}^{(0)} = \dfrac{ 2\varepsilon_{ \omega-1} }{ \varepsilon_{\omega+1}^2 - \gamma^2 \varepsilon_{ \omega-1}^2} r_i^2.
\end{align*}

We next compute the leading-order second-harmonic field \( u_{2\omega}^{(0)} \). Substituting the explicit expression \eqref{eq:u-omega-leading-explicit} into \eqref{eq:Ps-perp-i-leading},  \eqref{eq:Ps-perp-e-leading}, \eqref{eq:Ps-parallel-i-leading} and \eqref{eq:Ps-parallel-e-leading}, and then applying the surface-charge relations \eqref{eq:sigma-s-i-expansion} and \eqref{eq:sigma-s-e-expansion}, we obtain the following explicit expressions for the leading-order normal surface nonlinear polarizations and the associated surface charges
\begin{align}
    \label{eq:Ps-perp-i-leading-explicit} P_{\perp,i}^{s,(0)} & = E^2 \frac{ 8\chi_\perp^s \varepsilon_\omega^2 }{ (\varepsilon_{\omega+1}^2 - \gamma^2 \varepsilon_{\omega-1}^2)^2} \bigl( 1 + \cos( 2\theta ) \bigr) \triangleq E^2 \mathscr{A} \bigl( 1 + \cos( 2\theta ) \bigr), \\[4pt]
    \label{eq:Ps-perp-e-leading-explicit} P_{\perp,e}^{s,(0)} & = E^2 \frac{ 2\chi_\perp^s \varepsilon_\omega^2(\varepsilon_{\omega+1} - \gamma^2 \varepsilon_{\omega-1})^2}{ (\varepsilon_{ \omega+1}^2 -\gamma^2\varepsilon_{\omega-1}^2)^2} \bigl( 1 + \cos( 2\theta ) \bigr) \triangleq E^2 \mathscr{B} \bigl( 1 + \cos( 2\theta ) \bigr), \\[4pt]
    \label{eq:sigma-i-leading-explicit}
    \sigma_i^{s,(0)} & = E^2 \frac{ 32\chi_\parallel^s \varepsilon_\omega^2}{(\varepsilon_{\omega+1}^2 - \gamma^2 \varepsilon_{\omega-1}^2 )^2 } r_i^{-1} \cos( 2\theta ) \triangleq E^2 \mathscr{C} r_i^{-1} \cos( 2\theta ), \\[4pt]
    \label{eq:sigma-e-leading-explicit}
    \sigma_e^{s,(0)} & = E^2 \frac{ 8\chi_\parallel^s \varepsilon_\omega (\varepsilon_{\omega+1}^2 - \gamma^4 \varepsilon_{\omega-1}^2) }{ ( \varepsilon_{ \omega+1 }^2 - \gamma^2 \varepsilon_{\omega-1}^2 )^2} r_e^{-1} \cos( 2\theta) \triangleq E^2 \mathscr{D} r_e^{-1} \cos( 2\theta) .
\end{align}
Therefore, substituting \eqref{eq:Ps-perp-i-leading-explicit}--\eqref{eq:sigma-e-leading-explicit} into \eqref{eq:u-2omega-leading-density} and applying the explicit formulas \eqref{eq:single-layer-derivative-explicit}, \eqref{eq:double-layer-derivative-explicit} and \eqref{eq:NP-operator-explicit}, we obtain
\begin{align*}
    %\label{eq:u-2omega-leading-density-explicit}
    \begin{cases}
    \phi_i^{(0)} = 4\pi E^2 \mathscr{A} \bigl( 1 + \cos( 2\theta ) \bigr) & \text{on } \partial D, \\[8pt]
    \phi_e^{(0)} = -4\pi E^2 \mathscr{B} \bigl( 1 + \cos( 2\theta ) \bigr) & \text{on } \partial \Omega, \\[8pt]
    \psi_i^{(0)} = 4\pi E^2 r_i^{-1} \dfrac{ 2\gamma^2 \varepsilon_{ 2\omega-1}^2 \bigl( \gamma^2 \mathscr{A} -\mathscr{B} - \varepsilon_{2\omega-1}^{-1} \mathscr{D} \bigr) - 2\varepsilon_{2\omega+1} \varepsilon_{2\omega-1} \bigl( \mathscr{A} - \gamma^2 \mathscr{B} - \varepsilon_{2\omega-1}^{-1} \mathscr{C} \bigr) }{ \varepsilon_{2\omega+1}^2 - \gamma^4 \varepsilon_{2\omega-1}^2} \cos( 2\theta ) & \text{on } \partial D, \\[8pt]
    \psi_e^{(0)} = -4\pi E^2 r_e^{-1} \dfrac{2 \gamma^2 \varepsilon_{2\omega-1}^2 \bigl( \mathscr{A} - \gamma^2 \mathscr{B} - \varepsilon_{2\omega-1}^{-1} \mathscr{C} \bigr) - 2 \varepsilon_{2\omega+1} \varepsilon_{2\omega-1} \bigl( \gamma^2 \mathscr{A} - \mathscr{B} - \varepsilon_{2\omega-1}^{-1} \mathscr{D} \bigr) }{ \varepsilon_{2\omega+1}^2 - \gamma^4 \varepsilon_{2\omega-1}^2} \cos( 2\theta ) & \text{on } \partial \Omega.
    \end{cases}
\end{align*}
Substituting the densities $(\phi_i^{(0)},\phi_e^{(0)},\psi_i^{(0)},\psi_e^{(0)})$ into \eqref{eq:u-2omega-leading} and applying the explicit formulas \eqref{eq:single-layer-explicit}, \eqref{eq:double-layer-explicit} and \eqref{eq:double-layer-constant-density-explicit}, the leading-order solution for the second-harmonic field in the hollow nanowire with a perfect outer boundary is obtained after straightforward calculation as follows
\scriptsize
\begin{align}
    \label{eq:u-2omega-leading-explicit}
    \begin{aligned}
    u_{2\omega}^{(0)} = \begin{cases}
    \ds -4\pi E^2 ( \mathscr{B} - \mathscr{A} ) -4\pi E^2 \dfrac{ 2\varepsilon_{2\omega} ( \gamma^4 \varepsilon_{2\omega-1} - \varepsilon_{2\omega+1} ) \mathscr{A} + 4\varepsilon_{2\omega} \gamma^2 \mathscr{B} + ( \gamma^4 \varepsilon_{2\omega-1} + \varepsilon_{2\omega+1} ) \mathscr{C} -2\varepsilon_{2\omega} \gamma^2 \mathscr{D} }{ 2\left( \varepsilon_{2\omega+1}^2 - \gamma^4 \varepsilon_{2\omega-1}^2 \right) } \dfrac{r^2}{r_i^2}  \cos( 2\theta ), & r<r_i, \vspace{1em}\\
    \ds -4\pi E^2 \mathscr{B} - 4\pi E^2 \dfrac{ 2\varepsilon_{2\omega-1} \gamma^2 \mathscr{A} +2\varepsilon_{2\omega+1} \mathscr{B} +\varepsilon_{2\omega-1} \gamma^2 \mathscr{C} -\varepsilon_{2\omega+1} \mathscr{D} }{ 2\left( \varepsilon_{2\omega+1}^2 - \gamma^4 \varepsilon_{2\omega-1}^2 \right) } \dfrac{r^2}{r_e^2}  \cos( 2\theta ) & \vspace{1em}\\
    \quad \quad \quad \quad \ -4\pi E^2 \dfrac{ 2\varepsilon_{2\omega+1} \mathscr{A} + 2\varepsilon_{2\omega-1} \gamma^2 \mathscr{B} + \varepsilon_{2\omega+1} \mathscr{C} - \varepsilon_{2\omega-1} \gamma^2 \mathscr{D} }{ 2\left( \varepsilon_{2\omega+1}^2 - \gamma^4 \varepsilon_{2\omega-1}^2\right) } \dfrac{r_i^2}{r^2}  \cos(2\theta), & r_i<r<r_e, \vspace{1em} \\
    \ds -4\pi E^2 \dfrac{ 4\varepsilon_{2\omega} \gamma^2 \mathscr{A} +2\varepsilon_{2\omega} \left(\gamma^4  \varepsilon_{2\omega-1} - \varepsilon_{2\omega+1} \right) \mathscr{B} + 2\varepsilon_{2\omega} \gamma^2 \mathscr{C} - \left( \gamma^4 \varepsilon_{2\omega-1} + \varepsilon_{2\omega+1} \right) \mathscr{D} }{ 2\left( \varepsilon_{2\omega+1}^2 - \gamma^4 \varepsilon_{2\omega-1}^2 \right)} \dfrac{r_e^2}{r^2}  \cos( 2\theta ), & r>r_e.
    \end{cases}
    \end{aligned}
\end{align}
\normalsize

The above expression implies that, for a perfect hollow nanowire, the leading-order second-harmonic field in the exterior region is purely quadrupolar, containing only the \(\cos(2\theta)\) mode and decaying as \(O(r^{-2})\). More precisely,
\begin{align*}
    u_{2\omega}^{(0)}(x) = \hat{u}( \varepsilon_\omega, \varepsilon_{2\omega}, \gamma, \chi_\perp^s, \chi_\parallel^s) \dfrac{r_e^2}{r^2} \cos( 2\theta), \quad r > r_e,
\end{align*}
where \(C\) depends on the material parameters \( \varepsilon_\omega \), \( \varepsilon_{2\omega} \), the geometric ratio \( \gamma \), and the surface nonlinear susceptibilities \( \chi_\perp^s \) and \( \chi_\parallel^s \). In particular, no dipole term of order \( O( |x|^{-1} ) \) appears, and hence dipole radiation is forbidden.

This completes the proof.
\end{proof}

\begin{proof}[Proof of Corollary \ref{cor:inner-radius}]
From the exterior expression obtained in the proof of Theorem~\ref{thm:mainthm1}, the amplitude can be written as
\begin{align*}
 \hat{u}(\gamma)
 =-2\pi E^2\frac{(1-\gamma^2)\mathcal T(\gamma^2)}{\left(\varepsilon_{\omega+1}^2-\gamma^2\varepsilon_{\omega-1}^2\right)^2\left(\varepsilon_{2\omega+1}^2-\gamma^4\varepsilon_{2\omega-1}^2\right)}.
\end{align*}
Setting $\gamma=0$ gives the expression for $\hat{u}(0)$. Since $\hat{u}$ depends on $\gamma$ only through $\gamma^2$, a Taylor expansion at $\gamma=0$, gives
\begin{align*}
    \hat{u}(\gamma)=\hat{u}(0)+\mathcal O(\gamma^2),\qquad\gamma\to0.
\end{align*}
As $\gamma\to1^-$, the denominator factors satisfy
\begin{align*}
    \varepsilon_{\omega+1}^2 - \gamma^2 \varepsilon_{\omega-1}^2 \to 4\varepsilon_\omega, \quad \varepsilon_{2\omega+1}^2 - \gamma^4 \varepsilon_{2\omega-1}^2 \to 4 \varepsilon_{2\omega}
\end{align*}
Hence, provided $\varepsilon_\omega \varepsilon_{2\omega} \neq 0$, expansion about $\gamma=1$ yields
\begin{align*}
    \hat{u}( \gamma ) = \hat{u}_*(1-\gamma^2) + \mathcal O \bigl( (1-\gamma^2)^2 \bigr).
\end{align*}
Since \( 1-\gamma^2 = \mathcal O(r_e-r_i), \) the third assertion follows. Finally, for $0<\gamma<1$, $\hat{u}(\gamma)$ is smooth unless
\begin{align*}
    \varepsilon_{\omega+1}^2 - \gamma^2 \varepsilon_{\omega-1}^2 \to 0 \qquad \text{or} \qquad \varepsilon_{2\omega+1}^2 - \gamma^4 \varepsilon_{2\omega-1}^2 \to 0,
\end{align*}
which are precisely the stated plasmonic resonance conditions.

This completes the proof.
\end{proof}

We now proceed to the proof of Theorem \ref{thm:mainthm2}.

\begin{proof}[Proof of Theorem \ref{thm:mainthm2}]
Using \eqref{eq:boundary-term}, \eqref{eq:f-series} and \eqref{eq:u-omega-leading-explicit}, we obtain
\begin{align}
    \label{eq:I1}
    I_1 = {} & -E  \sum_{p\in\{-1, 1\}} r_e \mathscr{G} \cos \bigl( (n+p) \theta \bigl), \\
    \label{eq:I2}
    I_2 = {} & -E  \sum_{p\in\{-1, 1\}} \mathscr{H} p(n+p) \cos \bigl( (n+p) \theta \bigl),
\end{align}
where
\begin{align*}
    \mathscr{G} = \dfrac{ \gamma^2 \varepsilon_{\omega-1}^2 -
    \varepsilon_{\omega+1} \varepsilon_{\omega-1} }
    {\varepsilon_{\omega+1}^2 - \gamma^2 \varepsilon_{\omega-1}^2}, \quad \mathscr{H} = \dfrac{ \gamma^2 \varepsilon_{\omega-1}^2 +
    \varepsilon_{\omega+1} \varepsilon_{\omega-1} }
    {\varepsilon_{\omega+1}^2 - \gamma^2 \varepsilon_{\omega-1}^2}.
\end{align*}
Solving equation \eqref{eq:u-omega-first-density} with \eqref{eq:I1} and \eqref{eq:I2} yields
\begin{align}
    \label{eq:u-omega-first-density-explicit}
    \begin{cases}
    \vartheta_e^{(1)} = E \sum\limits_{p\in\{-1, 1\}} r_e \mathscr{G} \cos \bigl( (n+p) \theta \bigl) & \text{on } \partial \Omega, \\[12pt]
    \varphi_i^{(1)} = -\sum\limits_{p\in\{-1, 1\}} (n+p) E \dfrac{ ( 2\lambda_\omega-1) \gamma^{n+p-1} \mathscr{G} + 2p \gamma^{n+p-1} \varepsilon_{\omega-1}^{-1}  \mathscr{H}}{ 4\lambda_\omega^{2} - \gamma^{2(n+p)} } \cos((n+p) \theta)  & \text{on } \partial D, \\[12pt]
    \varphi_e^{(1)} = \sum\limits_{p\in\{-1, 1\}} (n+p) E \dfrac{ (2 \lambda_\omega - \gamma^{2(n+p)} ) \mathscr{G} + 4p \lambda_\omega \varepsilon_{\omega-1}^{-1} \mathscr{H} }{ 4\lambda_\omega^{2} - \gamma^{2(n+p)}} \cos((n+p)\theta)  & \text{on } \partial\Omega.
\end{cases}
\end{align}
Substituting \eqref{eq:u-omega-first-density-explicit} into \eqref{eq:u-omega-first} and applying the explicit formulas \eqref{eq:single-layer-explicit} and \eqref{eq:double-layer-explicit} yields the first-order solution for the fundamental field
\begin{align}
    \label{eq:u-omega-first-explicit}
    u_\omega^{(1)}(x) = \begin{cases}
    \ds E r_e \sum_{p\in\{-1,1\}} \alpha_i^{(n+p)} \left( \dfrac{r}{r_i} \right)^{n+p} \cos\bigl((n+p)\theta\bigr), & r<r_i, \\[10pt]
    \ds E r_e \sum_{p\in\{-1,1\}} \left[ \alpha_{m-}^{(n+p)} \left( \dfrac{r}{r_e} \right)^{n+p} + \alpha_{m+}^{(n+p)} \left( \dfrac{r_i}{r} \right)^{n+p} \right] \cos\bigl((n+p)\theta\bigr), & r_i<r<r_e, \\[10pt]
    \ds E r_e \sum_{p\in\{-1,1\}} \alpha_e^{(n+p)} \left( \dfrac{r_e}{r} \right)^{n+p} \cos\bigl((n+p)\theta\bigr), & r>r_e.
    \end{cases}
\end{align}
where the eight coefficients obtained from the layer-potential calculation are given explicitly by
\scriptsize
\begin{align*}
    \alpha_i^{(n+1)} & = \dfrac{4\gamma^{n+3}\varepsilon_\omega\varepsilon_{\omega-1}^{2}}{\left(\varepsilon_{\omega+1}^{2}-\gamma^{2}\varepsilon_{\omega-1}^{2}\right)\left(\varepsilon_{\omega+1}^{2}-\gamma^{2n+2}\varepsilon_{\omega-1}^{2}\right)}, &\qquad \alpha_i^{(n-1)} & = -\dfrac{4\gamma^{n-1}\varepsilon_\omega\varepsilon_{\omega+1}\varepsilon_{\omega-1}}{\left(\varepsilon_{\omega+1}^{2}-\gamma^{2}\varepsilon_{\omega-1}^{2}\right)\left(\varepsilon_{\omega+1}^{2}-\gamma^{2n-2}\varepsilon_{\omega-1}^{2}\right)}, \\[10pt]
    \alpha_{m-}^{(n+1)} & = \dfrac{2\gamma^{2}\varepsilon_{\omega+1}\varepsilon_{\omega-1}^{2}}{\left(\varepsilon_{\omega+1}^{2}-\gamma^{2}\varepsilon_{\omega-1}^{2}\right)\left(\varepsilon_{\omega+1}^{2}-\gamma^{2n+2}\varepsilon_{\omega-1}^{2}\right)}, &\qquad \alpha_{m+}^{(n+1)} & = \dfrac{2\gamma^{n+3}\varepsilon_{\omega-1}^{3}}{\left(\varepsilon_{\omega+1}^{2}-\gamma^{2}\varepsilon_{\omega-1}^{2}\right)\left(\varepsilon_{\omega+1}^{2}-\gamma^{2n+2}\varepsilon_{\omega-1}^{2}\right)}, \\[10pt]
    \alpha_{m-}^{(n-1)} & = -\dfrac{2\varepsilon_{\omega+1}^{2}\varepsilon_{\omega-1}}{\left(\varepsilon_{\omega+1}^{2}-\gamma^{2}\varepsilon_{\omega-1}^{2}\right)\left(\varepsilon_{\omega+1}^{2}-\gamma^{2n-2}\varepsilon_{\omega-1}^{2}\right)}, &\qquad \alpha_{m+}^{(n-1)} & = -\dfrac{2\gamma^{n-1}\varepsilon_{\omega+1}\varepsilon_{\omega-1}^{2}}{\left(\varepsilon_{\omega+1}^{2}-\gamma^{2}\varepsilon_{\omega-1}^{2}\right)\left(\varepsilon_{\omega+1}^{2}-\gamma^{2n-2}\varepsilon_{\omega-1}^{2}\right)}, \\[10pt]
    \alpha_e^{(n+1)} & = \dfrac{\varepsilon_{\omega+1}\varepsilon_{\omega-1}\left[\varepsilon_{\omega+1}^{2}-\left(\gamma^{2}+\gamma^{2n+2}-\gamma^{2n+4}\right)\varepsilon_{\omega-1}^{2}\right]}{\left(\varepsilon_{\omega+1}^{2}-\gamma^{2}\varepsilon_{\omega-1}^{2}\right)\left(\varepsilon_{\omega+1}^{2}-\gamma^{2n+2}\varepsilon_{\omega-1}^{2}\right)}, &\qquad \alpha_e^{(n-1)} & = \dfrac{\varepsilon_{\omega-1}^{2}\left[\gamma^{2n}\varepsilon_{\omega-1}^{2}+\left(1-\gamma^{2}-\gamma^{2n-2}\right)\varepsilon_{\omega+1}^{2}\right]}{\left(\varepsilon_{\omega+1}^{2}-\gamma^{2}\varepsilon_{\omega-1}^{2}\right)\left(\varepsilon_{\omega+1}^{2}-\gamma^{2n-2}\varepsilon_{\omega-1}^{2}\right)}.
\end{align*}
\normalsize

We now proceed to compute the first-order solution \(u_{2\omega}^{(1)}\) of the second-harmonic field. To this end, we first derive the explicit expressions for the source terms \(I_3\) and \(I_4\), which are required in the boundary conditions. We first compute the source term~$I_3$. Recall that
\begin{align*}
    I_3 = f\left( \frac{\partial u_{2\omega}^{(0)}}{\partial \nu_e} \Big|_{-} - \frac{\partial u_{2\omega}^{(0)}}{\partial \nu_e} \Big|_{+} \right) + 4\pi P_{\perp,e}^{s,(1)}.
\end{align*}
For the first term, it follows from \eqref{eq:f-series} and \eqref{eq:u-2omega-leading-explicit} that
\begin{align*}
    f \left(  \frac{\partial u_{2\omega}^{(0)} }{\partial \nu_e} \Big|_{-} - \frac{\partial u_{2\omega}^{(0)} }{ \partial \nu_e } \Big|_{+} \right) & = -4\pi E^2 \frac{ 2\gamma^2 \varepsilon_{2\omega-1} \mathscr{A} + \bigl( \gamma^4 \varepsilon_{2\omega-1} - \varepsilon_{2\omega+1} \bigr) \varepsilon_{2\omega-1} \mathscr{B} + \gamma^2 \varepsilon_{2\omega-1} \mathscr{C} - \varepsilon_{2\omega+1} \mathscr{D} }{ \varepsilon_{2\omega+1}^2 - \gamma^4\varepsilon_{2\omega-1}^2 } \\
    &\qquad \quad \times \left[ \cos \bigl( (n-2) \theta \bigr) + \cos \bigl( (n+2) \theta \bigr) \right].
\end{align*}
For the surface nonlinear polarization term, and using \eqref{eq:Ps-perp-e-first}, \eqref{eq:f-series}, \eqref{eq:u-omega-leading-explicit} and \eqref{eq:u-omega-first-explicit}, we obtain
\scriptsize
\begin{align*}
    4\pi P_{\perp,e}^{s,(1)} = {} & -2\pi E^2 \chi_\perp^s \varepsilon_\omega^2 \alpha_e^- \left[ 2(n-1) \left( \alpha_{m-}^{(n-1)} - \gamma^{n-1} \alpha_{m+}^{(n-1)} \right) + n\alpha_{m-}^{(0)} + (n-2) r_e^{-2} \alpha_{m+}^{(0)} \right] \cos \bigl( (n-2) \theta \bigr) \\[4pt]
    & - 4\pi E^2 \chi_\perp^s \varepsilon_\omega^2 \alpha_e^- \left[ (n-1) \left( \alpha_{m-}^{(n-1)} - \gamma^{n-1} \alpha_{m+}^{(n-1)} \right) + (n+1) \left( \alpha_{m-}^{(n+1)} - \gamma^{n+1} \alpha_{m+}^{(n+1)} \right) - 2r_e^{-2} \alpha_{m+}^{(0)} \right] \cos( n\theta ) \\[4pt]
    & - 2\pi E^2 \chi_\perp^s \varepsilon_\omega^2 \alpha_e^- \left[ 2(n+1) \left( \alpha_{m-}^{(n+1)} - \gamma^{n+1} \alpha_{m+}^{(n+1)} \right) - n\alpha_{m-}^{(0)} - (n+2) r_e^{-2} \alpha_{m+}^{(0)} \right] \cos \bigl((n+2)\theta\bigr).
\end{align*}
\normalsize
Here,
\begin{equation*}
    \alpha_e^{\pm} = 2 \frac{\varepsilon_{\omega+1} \pm \gamma^2 \varepsilon_{\omega-1}}{\varepsilon_{\omega+1}^2 - \gamma^2 \varepsilon_{\omega-1}^2}.
\end{equation*}
Consequently,
\begin{equation}
    \label{eq:I3}
    I_3 = \sum_{p\in\{-2,0,2\}} \mathcal{I}_3^{(n+p)} \cos\bigl((n+p)\theta\bigr),
\end{equation}
where
\begin{align*}
    \mathcal{I}_3^{(n-2)} = {} & -4\pi E^2 \frac{ 2\gamma^2 \varepsilon_{ 2\omega - 1} \mathscr{A} + \left( \gamma^4 \varepsilon_{2\omega - 1} - \varepsilon_{2\omega + 1} \right) \varepsilon_{ 2\omega - 1} \mathscr{B} + \gamma^2 \varepsilon_{ 2\omega - 1} \mathscr{C} - \varepsilon_{ 2\omega + 1} \mathscr{D} }{ \varepsilon_{ 2\omega + 1}^2 -
    \gamma^4 \varepsilon_{ 2\omega - 1}^2 } \\[4pt]
    & \quad - 2\pi E^2 \chi_\perp^s \varepsilon_\omega^2 \alpha_e^- \left[ 2(n-1) \left( \alpha_{m-}^{(n-1)} - \gamma^{n-1} \alpha_{m+}^{(n-1)} \right) + n \alpha_{m-}^{(0)} + (n-2) r_e^{-2} \alpha_{m+}^{(0)} \right], \\[8pt]
    \mathcal{I}_3^{(n)} ={}& -4\pi E^2\chi_\perp^s\varepsilon_\omega^2 \alpha_e^-\left[ (n-1) \left( \alpha_{m-}^{(n-1)} - \gamma^{n-1}\alpha_{m+}^{(n-1)} \right) + (n+1) \left( \alpha_{m-}^{(n+1)} - \gamma^{n+1} \alpha_{m+}^{(n+1)} \right) - 2r_e^{-2}\alpha_{m+}^{(0)} \right],\\[4pt]
    \mathcal{I}_3^{(n+2)} ={}& -4\pi E^2 \frac{ 2\gamma^2\varepsilon_{2\omega-1}\mathscr{A} + \bigl( \gamma^4\varepsilon_{2\omega-1} - \varepsilon_{2\omega+1} \bigr) \varepsilon_{2\omega-1}\mathscr{B} + \gamma^2\varepsilon_{2\omega-1}\mathscr{C} - \varepsilon_{2\omega+1}\mathscr{D} }{ \varepsilon_{2\omega+1}^2 - \gamma^4\varepsilon_{2\omega-1}^2 } \\[4pt]
    & \quad -2\pi E^2 \chi_\perp^s \varepsilon_\omega^2
    \alpha_e^- \left[ 2(n+1) \left( \alpha_{m-}^{(n+1)} - \gamma^{n+1} \alpha_{m+}^{(n+1)} \right) - n\alpha_{m-}^{(0)} - (n+2) r_e^{-2} \alpha_{m+}^{(0)} \right].
\end{align*}

We next compute the source term~$I_4$. Recall that
\begin{align*}
    I_4 = {} f \left( \varepsilon_{2\omega} \frac{ \partial^2 u_{2\omega}^{(0)}}{\partial\nu_e^2} \Big|_- -  \frac{\partial^2 u_{2\omega }^{(0)}}{\partial\nu_e^2} \Big|_+ \right) + f' \left( \frac{\partial u_{2\omega }^{(0)}}{\partial\tau_e} \Big|_+ - \varepsilon_{ 2\omega} \frac{\partial u_{2\omega }^{(0)}}{ \partial \tau_e} \Big|_- \right) - 4\pi \sigma_e^{ s,(1)}.
\end{align*}
For the first term, direct calculation gives
\begin{align*}
    f \left( \varepsilon_{2\omega}  \frac{\partial^2 u_{2\omega}^{(0)}}{\partial\nu_e^2} \Big|_- - \frac{\partial^2 u_{2\omega}^{(0)} }{\partial \nu_e^2} \Big|_+ \right) = & -\frac{2\pi E^2}{r_e} \frac{1} { \varepsilon_{2\omega+1}^2 - \gamma^4 \varepsilon_{2\omega-1}^2 } \Big[ 8 \gamma^2 \varepsilon_{2\omega} \varepsilon_{2\omega-1} \mathscr{A} + 8 \varepsilon_{2\omega} \varepsilon_{2\omega+1} \mathscr{B} + 4 \gamma^2 \varepsilon_{2\omega} \varepsilon_{2\omega-1} \mathscr{C} \\[4pt]
    & + \left( \varepsilon_{2\omega+1} \bigl( 3 - \varepsilon_{2\omega} \bigr) - 3 \gamma^4 \varepsilon_{2\omega-1}^2 \right) \mathscr{D} \Big] \times \left[ \cos\bigl((n-2)\theta\bigr) + \cos\bigl((n+2)\theta\bigr) \right].
\end{align*}
For the second term, we obtain
\begin{align*}
    f' \left( \frac{\partial u_{2\omega}^{(0)} }{ \partial\tau_e} \Big|_+ - \varepsilon_{ 2\omega } \frac{\partial  u_{2\omega}^{(0)} }{ \partial \tau_e} \Big|_- \right) ={} & \frac{2\pi E^2 n}{r_e} \frac{1}{ \varepsilon_{ 2\omega+1}^2 - \gamma^4 \varepsilon_{2\omega-1}^2 } \Big[ 4\gamma^2 \varepsilon_{2\omega} \varepsilon_{ 2\omega-1} \mathscr{A} + 4\varepsilon_{2\omega} \varepsilon_{2\omega+1} \mathscr{B} + 2\gamma^2 \varepsilon_{2\omega} \varepsilon_{2\omega-1} \mathscr{C} \\[4pt]
    & - \varepsilon_{2\omega-1} \left( \varepsilon_{2\omega+1} + \gamma^4 \varepsilon_{2\omega-1} \right) \mathscr{D} \Big] \times \left[ \cos\bigl((n-2)\theta\bigr) - \cos\bigl((n+2)\theta\bigr) \right].
\end{align*}
It remains to compute the first-order surface charge $\sigma_e^{s,(1)}$. Substituting the explicit expressions for $u_\omega^{(0)}$ and $u_\omega^{(1)}$ into \eqref{eq:Ps-parallel-e-leading}--\eqref{eq:sigma-s-e-expansion} yields
\begin{align*}
 -4\pi\sigma_e^{s,(1)}={}&\frac{4\pi E^2}{r_e}\chi_\parallel^s\varepsilon_\omega\Big[2(n-1)(n-2)\Big(\gamma^{n-1}\alpha_{m-}^{(0)}\alpha_{m+}^{(n-1)}-r_e^{-2}\alpha_{m+}^{(0)}\alpha_{m-}^{(n-1)}\Big)\\[4pt]
 &\hspace{18mm}-(n^2-2n-1)\bigl(\alpha_{m-}^{(0)}\bigr)^2-(n^2-4n+5)r_e^{-4}\bigl(\alpha_{m+}^{(0)}\bigr)^2\Big]\cos\bigl((n-2)\theta\bigr)\\[6pt]
 &+\frac{8\pi E^2n}{r_e}\chi_\parallel^s\varepsilon_\omega\Big[(n-1)\Big(\alpha_{m-}^{(0)}\alpha_{m-}^{(n-1)}-\gamma^{n-1}r_e^{-2}\alpha_{m+}^{(0)}\alpha_{m+}^{(n-1)}\Big)\\[4pt]
 &\hspace{18mm}+(n+1)\Big(\gamma^{n+1}\alpha_{m-}^{(0)}\alpha_{m+}^{(n+1)}-r_e^{-2}\alpha_{m+}^{(0)}\alpha_{m-}^{(n+1)}\Big)+2nr_e^{-2}\alpha_{m-}^{(0)}\alpha_{m+}^{(0)}\Big]\cos(n\theta)\\[6pt]
 &-\frac{4\pi E^2}{r_e}\chi_\parallel^s\varepsilon_\omega\Big[(n^2+2n-1)\bigl(\alpha_{m-}^{(0)}\bigr)^2+(n^2+4n+5)r_e^{-4}\bigl(\alpha_{m+}^{(0)}\bigr)^2\\[4pt]
 &\hspace{18mm}-2(n+1)(n+2)\Big(\alpha_{m-}^{(0)}\alpha_{m-}^{(n+1)}-\gamma^{n+1}r_e^{-2}\alpha_{m+}^{(0)}\alpha_{m+}^{(n+1)}\Big)\Big]\cos\bigl((n+2)\theta\bigr).
\end{align*}
Combining the above three terms, we obtain
\begin{align}
 \label{eq:I4}
 I_4=\frac1{r_e}\sum_{p\in\{-2,0,2\}}\mathcal{I}_4^{(n+p)}\cos\bigl((n+p)\theta\bigr),
\end{align}
where
\scriptsize
\begin{align*}
 \mathcal{I}_4^{(n-2)}={}&2\pi E^2(n-2) \frac{ 4\gamma^2 \varepsilon_{2\omega} \varepsilon_{2\omega-1} \mathscr A + 4\varepsilon_{2\omega} \varepsilon_{2\omega+1} \mathscr B + 2\gamma^2 \varepsilon_{2\omega} \varepsilon_{2\omega-1} \mathscr C - \varepsilon_{2\omega-1} ( \varepsilon_{2\omega+1} + \gamma^4 \varepsilon_{2\omega-1}) \mathscr D}{\varepsilon_{2\omega+1}^2 - \gamma^4 \varepsilon_{2\omega-1}^2} - 2\pi E^2\mathscr D\\[4pt]
 &+4\pi E^2 \chi_\parallel^s \varepsilon_\omega \Big[2(n-1)(n-2) \Big( \gamma^{n-1} \alpha_{m-}^{(0)} \alpha_{m+}^{(n-1)} - r_e^{-2} \alpha_{m+}^{(0)} \alpha_{m-}^{(n-1)} \Big)-(n^2-2n-1) \bigl(\alpha_{m-}^{(0)} \bigr)^2 - (n^2-4n+5) r_e^{-4} \bigl( \alpha_{m+}^{(0)} \bigr)^2\Big],\\[8pt]
 \mathcal{I}_4^{(n)}={}&8\pi E^2n \chi_\parallel^s \varepsilon_\omega \Big[(n-1) \Big(\alpha_{m-}^{(0)} \alpha_{m-}^{(n-1)} - \gamma^{n-1}r_e^{-2} \alpha_{m+}^{(0)} \alpha_{m+}^{(n-1)}\Big)\\[4pt]
 &\hspace{18mm}+(n+1)\Big(\gamma^{n+1} \alpha_{m-}^{(0)} \alpha_{m+}^{(n+1)} - r_e^{-2} \alpha_{m+}^{(0)} \alpha_{m-}^{(n+1)} \Big)+2nr_e^{-2} \alpha_{m-}^{(0)}\alpha_{m+}^{(0)}\Big],\\[8pt]
 \mathcal{I}_4^{(n+2)}={}&-2\pi E^2(n+2) \frac{4\gamma^2 \varepsilon_{2\omega} \varepsilon_{2\omega-1} \mathscr A + 4\varepsilon_{2\omega} \varepsilon_{2\omega+1} \mathscr B + 2\gamma^2 \varepsilon_{2\omega} \varepsilon_{2\omega-1} \mathscr C - \varepsilon_{2\omega-1} (\varepsilon_{2\omega+1} + \gamma^4 \varepsilon_{2\omega-1}) \mathscr D}{\varepsilon_{2\omega+1}^2 - \gamma^4 \varepsilon_{2\omega-1}^2} - 2\pi E^2 \mathscr D\\[4pt]
 &-4\pi E^2 \chi_\parallel^s \varepsilon_\omega \Big[(n^2+2n-1) \bigl( \alpha_{m-}^{(0)}\bigr)^2 + (n^2+4n+5)r_e^{-4}\bigl(\alpha_{m+}^{(0)}\bigr)^2 - 2(n+1)(n+2) \Big( \alpha_{m-}^{(0)} \alpha_{m-}^{(n+1)} - \gamma^{n+1}r_e^{-2} \alpha_{m+}^{(0)} \alpha_{m+}^{(n+1)}\Big)\Big].
\end{align*}
\normalsize

To complete the computation of \(u_{2\omega}^{(1)}\), we next determine \(P_{\perp,i}^{s,(1)}\) and \(\sigma_i^{s,(1)}\) on the inner boundary \(\partial D\). Applying \eqref{eq:Ps-perp-i-first}, \eqref{eq:Ps-parallel-i-first} and \eqref{eq:sigma-s-i-expansion}, we obtain
\begin{align}
    \label{eq:Ps-perp-i-first-explicit}
    -4\pi P_{\perp,i}^{s,(1)} & = \sum_{p\in\{-2,0,2\}}\mathcal{M}^{(n+p)}\cos\bigl((n+p)\theta\bigr),\\[4pt]
    \label{eq:sigma-i-first-explicit}
    4\pi\sigma_i^{s,(1)} & = \sum_{p\in\{0,2\}}\frac1{r_i}\mathcal{N}^{(n+p)}\cos\bigl((n+p)\theta\bigr),
\end{align}
where
\begin{align*}
    \mathcal{M}^{(n-2)} & = 16\pi E^2\frac{\chi_\perp^s \varepsilon_\omega^2(n-1)}{\varepsilon_{\omega+1}^2 - \gamma^2 \varepsilon_{\omega-1}^2} \left( \gamma^{n-2} \alpha_{m-}^{(n-1)} - \gamma^{-1} \alpha_{m+}^{(n-1)} \right),\\[6pt]
    \mathcal{M}^{(n+2)} & = 16\pi E^2\frac{\chi_\perp^s \varepsilon_\omega^2(n+1)}{\varepsilon_{\omega+1}^2 - \gamma^2 \varepsilon_{\omega-1}^2} \left( \gamma^{n}\alpha_{m-}^{(n+1)} - \gamma^{-1} \alpha_{m+}^{(n+1)} \right),\\[6pt]
    \mathcal{M}^{(n)} & = \mathcal{M}^{(n-2)} + \mathcal{M}^{(n+2)},\\[6pt]
    \mathcal{N}^{(n)} & = -32\pi E^2\frac{\chi_\parallel^s \varepsilon_\omega^2 n(n-1)}{\varepsilon_{\omega+1}^2 -\gamma^2 \varepsilon_{\omega-1}^2} \left( \gamma^{n-2} \alpha_{m-}^{(n-1)} -\gamma^{-1} \alpha_{m+}^{(n-1)} \right),\\[6pt]
    \mathcal{N}^{(n+2)} & = -32\pi E^2\frac{\chi_\parallel^s \varepsilon_\omega^2 (n+1)(n+2)}{\varepsilon_{\omega+1}^2 - \gamma^2 \varepsilon_{\omega-1}^2} \left( \gamma^{n} \alpha_{m-}^{(n+1)} -\gamma^{-1} \alpha_{m+}^{(n+1)}\right).
\end{align*}

With the outer-boundary source terms \(I_3\) and \(I_4\) given by \eqref{eq:I3} and \eqref{eq:I4}, respectively, and the inner-boundary quantities \(P_{\perp,i}^{s,(1)}\) and \(\sigma_i^{s,(1)}\) given by \eqref{eq:Ps-perp-i-first-explicit} and \eqref{eq:sigma-i-first-explicit}, respectively, we now solve the coupled boundary integral system \eqref{eq:u-2omega-first-density} for the density quadruple \((\phi_i^{(1)}, \phi_e^{(1)}, \psi_i^{(1)}, \psi_e^{(1)})\). A direct calculation yields
\begin{align}
    \label{eq:u-2omega-first-density-explicit}
    \begin{cases}
    \ds \phi_i^{(1)} = -\sum\limits_{p\in\{-2,0,2\}} \mathcal{M}^{(n+p)} \cos \bigl((n+p)\theta\bigr) & \text{on } \partial D, \ms \\
    \ds \phi_e^{(1)} = -\sum\limits_{p\in\{-2,0,2\}} \mathcal{I}_3^{(n+p)} \cos \bigl((n+p)\theta\bigr) & \text{on } \partial \Omega, \ms \\
    \ds \psi_i^{(1)} = \sum\limits_{p\in\{-2,0,2\}} \dfrac{ 2r_e^{-1} \varepsilon_{2\omega-1} \left( \gamma^{-1} \varepsilon_{2\omega+1} \mathscr{P}^{(n+p)} + \gamma^{n+p-1} \varepsilon_{2\omega-1} \mathscr{Q}^{(n+p)} \right) }{ \varepsilon_{2\omega+1}^2 - \gamma^{2(n+p)} \varepsilon_{2\omega-1}^2 } \cos \bigl( (n+p) \theta\bigr) & \text{on } \partial D, \ms \\
    \ds \psi_e^{(1)} = \sum\limits_{p\in\{-2,0,2\}} \dfrac{ 2r_e^{-1}\varepsilon_{2\omega-1} \left( \varepsilon_{2\omega+1} \mathscr{Q}^{(n+p)} + \gamma^{n+p} \varepsilon_{2\omega-1} \mathscr{P}^{(n+p)} \right) }{ \varepsilon_{2\omega+1}^2 - \gamma^{2(n+p)} \varepsilon_{2\omega-1}^2 } \cos \bigl( (n+p) \theta \bigr) & \text{on } \partial \Omega,
    \end{cases}
\end{align}
where, for $p\in\{-2,0,2\}$,
\begin{align*}
    \mathscr{P}^{(n+p)} & := \frac{n+p}{2} \mathcal{M}^{(n+p)} + \frac{n+p}{2} \gamma^{n+p} \mathcal{I}_3^{(n+p)} + \frac{1} {\varepsilon_{2\omega-1}} \mathcal{N}^{(n+p)},\\[6pt]
    \mathscr{Q}^{(n+p)} & := -\frac{n+p}{2} \gamma^{n+p} \mathcal{M}^{(n+p)} - \frac{n+p}{2} \mathcal{I}_3^{(n+p)} + \frac{1}{\varepsilon_{2\omega-1}} \mathcal{I}_4^{(n+p)}.
\end{align*}
Here, we set $\mathcal{N}^{(n-2)}:=0$, since the corresponding mode is absent from the first-order inner surface charge.

Finally, substituting the density formulas \eqref{eq:u-2omega-first-density-explicit} into the layer-potential representation \eqref{eq:u-2omega-first} of \(u_{2\omega}^{(1)}\), and applying the explicit formulas for the single- and double-layer potentials on a disk established in Lemma~\ref{Layer-potentials-on-disk}, we obtain, for \(r>r_e\) and every integer \(n\geq3\),
\begin{align}
    \label{eq:u-2omega-first-explicit}
    u_{2\omega}^{(1)}(x)=\sum_{p\in\{-2,0,2\}}C^{(n+p)}\left(\frac{r_e}{r}\right)^{n+p}\cos\bigl((n+p)\theta\bigr),
\end{align}
where, for \(p\in\{-2,0,2\}\),
\scriptsize
\begin{align*}
     C^{(n+p)} = \frac{2\varepsilon_{2\omega} \gamma^{n+p} \mathcal{M}^{(n+p)} + \varepsilon_{2\omega} \left( \varepsilon_{2\omega+1} - \gamma^{2(n+p)} \varepsilon_{2\omega-1}\right) \mathcal{I}_3^{(n+p)} - \frac{ 2\varepsilon_{2\omega} \gamma^{n+p}}{n+p} \mathcal{N}^{(n+p)} - \frac{ \varepsilon_{2\omega+1} + \gamma^{2(n+p)} \varepsilon_{2\omega-1}}{n+p} \mathcal{I}_4^{(n+p)}}{\varepsilon_{2\omega+1}^2 - \gamma^{2(n+p)} \varepsilon_{2\omega-1}^2}.
\end{align*}
\normalsize
Equation~\eqref{eq:u-2omega-first-explicit} shows that the exterior first-order correction \(u_{2\omega}^{(1)}\) contains only angular modes of orders \(n-2\), \(n\), and \(n+2\), whose radial factors decay as \(r^{-(n-2)}\), \(r^{-n}\), and \(r^{-(n+2)}\), respectively. Provided that \(C^{(n-2)}\neq0\), the \((n-2)\)-nd mode determines the leading far-field behavior. It follows that:
\begin{enumerate}[label=\normalfont{(\roman*)}]
    \item If \(n=3\), the leading term decays as \(r^{-1}\), and the first-order correction exhibits dipole radiation.

    \item If \(n=4\), the leading term decays as \(r^{-2}\), and the first-order correction exhibits quadrupole radiation.

    \item If \(n\geq5\), the leading term decays as \(r^{-(n-2)}\), corresponding to a \(2^{\,n-2}\)-pole radiation field.
\end{enumerate}
Among the perturbations \(f(\theta)=r_e\cos(n\theta)\), \(n\geq3\), and under the above nonvanishing assumption, the first-order correction exhibits dipole radiation if and only if \(n=3\), corresponding to the lowest admissible symmetry degree \(2n=6\).

This completes the proof.

\end{proof}

Finally, we give the proof of Theorem~\ref{thm:mainthm3}.
\begin{proof}
Define
\begin{align*}
    \mathcal R_1&:= \frac1{ \left| \varepsilon_{\omega+1}^2 - \gamma^{2} \varepsilon_{\omega-1}^2 \right|^2 \left| \varepsilon_{2\omega+1}^2 - \gamma^{4} \varepsilon_{2\omega-1}^2 \right|},&
    \mathcal R_2&:= \frac1{ \left| \varepsilon_{\omega+1}^2 - \gamma^{2} \varepsilon_{\omega-1}^2 \right|^2},\\[6pt]
    \mathcal R_3&:= \frac1{ \left| \varepsilon_{\omega+1}^2 - \gamma^{2} \varepsilon_{\omega-1}^2 \right|^2 \left| \varepsilon_{\omega+1}^2 - \gamma^{2(n-1)} \varepsilon_{\omega-1}^2 \right|},&
    \mathcal R_4&:= \frac1{ \left| \varepsilon_{\omega+1}^2 - \gamma^{2} \varepsilon_{\omega-1}^2 \right|^2 \left| \varepsilon_{\omega+1}^2 - \gamma^{2(n+1)} \varepsilon_{\omega-1}^2 \right|}.
\end{align*}
The explicit source coefficients satisfy
\begin{align*}
    \mathcal I_3^{(n-2)},\ \mathcal I_4^{(n-2)}&= \mathcal O(\mathcal R_1 + \mathcal R_2 + \mathcal R_3),\\
    \mathcal I_3^{(n)},\ \mathcal I_4^{(n)}&= \mathcal O(\mathcal R_2 + \mathcal R_3 + \mathcal R_4),\\
    \mathcal I_3^{(n+2)},\ \mathcal I_4^{(n+2)}&= \mathcal O(\mathcal R_1 + \mathcal R_2 + \mathcal R_4).
\end{align*}
and
\begin{align*}
    \mathcal M^{(n-2)} &= \mathcal O(\mathcal R_3),&
    \mathcal M^{(n)} &= \mathcal O(\mathcal R_3 + \mathcal R_4),&
    \mathcal M^{(n+2)} &= \mathcal O(\mathcal R_4),\\
    \mathcal N^{(n-2)} &= 0,&
    \mathcal N^{(n)} &= \mathcal O(\mathcal R_3),&
    \mathcal N^{(n+2)} &= \mathcal O(\mathcal R_4).
\end{align*}
Substitution into the formula for $C^{(n+p)}$ following \eqref{eq:u-2omega-first-explicit} yields
\begin{align*}
    |C^{(n-2)}|& \lesssim\frac1{ \left| \varepsilon_{\omega+1}^2 - \gamma^{2} \varepsilon_{\omega-1}^2 \right|^2 \left| \varepsilon_{2\omega+1}^2 - \gamma^{2(n-2)} \varepsilon_{2\omega-1}^2 \right|} \left(1+\frac1{\left| \varepsilon_{2\omega+1}^2 - \gamma^{4} \varepsilon_{2\omega-1}^2 \right|} + \frac1{ \left| \varepsilon_{\omega+1}^2 - \gamma^{2(n-1)} \varepsilon_{\omega-1}^2 \right|} \right),\\[6pt]
    |C^{(n)}|& \lesssim\frac1{ \left| \varepsilon_{\omega+1}^2 - \gamma^{2} \varepsilon_{\omega-1}^2 \right|^2 \left| \varepsilon_{2\omega+1}^2 - \gamma^{2n} \varepsilon_{2\omega-1}^2 \right|} \left( 1 + \frac1{ \left| \varepsilon_{\omega+1}^2 - \gamma^{2(n-1)} \varepsilon_{\omega-1}^2\right|} + \frac1{\left| \varepsilon_{\omega+1}^2 - \gamma^{2(n+1)} \varepsilon_{\omega-1}^2 \right|} \right),\\[6pt]
    |C^{(n+2)}|& \lesssim\frac1{ \left| \varepsilon_{\omega+1}^2 - \gamma^{2} \varepsilon_{\omega-1}^2 \right|^2 \left| \varepsilon_{2\omega+1}^2 -\gamma^{2(n+2)} \varepsilon_{2\omega-1}^2 \right|} \left( 1 + \frac1{ \left|\varepsilon_{2\omega+1}^2 -\gamma^{4} \varepsilon_{2\omega-1}^2 \right|} + \frac1{ \left| \varepsilon_{\omega+1}^2 - \gamma^{2(n+1)} \varepsilon_{\omega-1}^2 \right|} \right).
\end{align*}
The common squared fundamental denominator gives the order-two resonance in the first case of \textup{(i)}, while the two higher fundamental denominators give order-one resonances. At the second-harmonic frequency, the source and lowest output denominators coincide precisely when $n=4$, yielding order $1+\delta_{n=4}$; the other output resonances have order one.

For \textup{(ii)}, the common fundamental factor multiplies every second-harmonic singular contribution. The higher fundamental factor with exponent $n-1$ occurs only in the lowest and middle modes, whereas that with exponent $n+1$ occurs only in the middle and highest modes. A higher fundamental factor therefore multiplies an output factor when both occur in the same mode; otherwise their contributions add. At the source resonance with $n=4$, the additional output factor multiplies both contributions in the lowest mode for the exponent $n-1$. For the exponent $n+1$, however, the repeated source factor and the fundamental factor occur in different modes and contribute additively.

This completes the proof.
\end{proof}

\section{Conclusion}
\label{sec:conclusion}
In this paper, we have developed a quantitative mathematical framework for second-harmonic generation in plasmonic hollow nanowires with a concentric circular cavity and a slightly perturbed outer boundary. By combining layer-potential methods with asymptotic expansions in the boundary perturbation parameter, we have derived explicit expressions for the fundamental and second-harmonic fields. Our analysis reveals the distinct roles of geometric symmetry breaking and the cavity-induced topological defect in second-harmonic generation. In particular, the outer-boundary symmetry determines the multipolar structure of the first-order second-harmonic response, while the cavity radius influences the radiation amplitudes and plasmonic resonance behavior. We have also characterized the effects of single and simultaneous resonances, demonstrating the role of the cavity geometry in controlling nonlinear optical enhancement. These results provide a quantitative understanding of the interplay between geometric symmetry, cavity structure, and nonlinear optical responses. Several directions remain open for further investigation. These include general boundary perturbations, nonconcentric or noncircular cavities, multilayer structures, and three-dimensional hollow nanoparticles. It would also be of interest to extend the present analysis to higher-order shape perturbations, higher-harmonic generation, and full-wave Maxwell models incorporating retardation effects. Such extensions may provide further insights into the geometric optimization and inverse design of nonlinear plasmonic nanostructures.

\end{document}